\documentclass[12pt,twoside]{article}
\usepackage{theorem}
\usepackage{fancyheadings}
\usepackage{latexsym,amssymb}

\newcommand{\eq}{\begin{equation}}
\newcommand{\en}{\end{equation}}
\newcommand{\re}[1]{\mbox{(\ref{#1})}}

\newtheorem{Theorem}{Theorem}
\newtheorem{theorem}[Theorem]{Theorem}
\newtheorem{lemma}[Theorem]{Lemma}

\newtheorem{corollary}[Theorem]{Corollary}
\newtheorem{construction}[Theorem]{Construction}
\newtheorem{proposition}[Theorem]{Proposition}
\newtheorem{example}[Theorem]{Example}
\newtheorem{exercise}[Theorem]{Exercise}
\newtheorem{defn}[Theorem]{Definition}

\newtheorem{question}[Theorem]{Question}
\newtheorem{conjecture}[Theorem]{Conjecture}
\newtheorem{condition}[Theorem]{Condition}
\newtheorem{remark}[Theorem]{Remark}
\newtheorem{problem}[Theorem]{Problem}
\newtheorem{jpfigure}[Theorem]{Figure}
\newtheorem{jptable}[Theorem]{Table}

\def\proof{\noindent{\bf Proof.\ \ }}

\def\endpf{\hfill $\Box$ \vskip .25in}

\newfont{\msbm}{msbm10 at 12pt}
\newfont{\eusb}{eusb10}
\newfont{\eusm}{eusm10}
\newfont{\eurb}{eurb10}
\newfont{\eurm}{eurm10}
\newfont{\eufb}{eufb10}
\newfont{\eufm}{eufm10}

\newcommand {\PR} {\mathbb{P}}
\newcommand {\ER} {\mathbb{E}}

\newcommand {\complex} {\mathbb{C}}

\newcommand {\IN} {\mbox{\msbm\symbol{'116}}}

\newcommand {\NN} {\mbox{\msbm\symbol{'116}}}

\newcommand{\te}{\rightarrow}
\newcommand{\ed}{\mbox{$ \ \stackrel{d}{=}$ }}

\newcommand{\eps}{\varepsilon}
\newcommand{\EE}{{\cal E}}

\newcommand{\Lev}{L\'evy}

\newcommand{\LK}{L\'evy-Khintchine}

\newcommand{\comment}[1]{}
\newcommand{\noshowcomment}{\renewcommand{\comment}[1]{}}

\newcommand{\lb}[1]{\label{#1}}

\newenvironment{thm}[1]{\begin{theorem}\label{#1}}{\end{theorem}}

\newenvironment{lmm}[1]{\begin{lemma}\label{#1}}{\end{lemma}}

\newenvironment{crl}[1]{\begin{corollary}\protect\label{#1}}{\end{corollary}}
\newenvironment{dfn}[1]{\begin{defn}\protect\label{#1}}{\end{defn}}

\newenvironment{prp}[1]{\begin{proposition}\protect\label{#1}}{\end{proposition}}

\noshowcomment

\newcommand{\keywords}[1]{\addvspace{\baselineskip}\noindent{{\bf Keywords.} #1}}
\newcommand{\affiliation}[1]{\addvspace{\baselineskip}\noindent\emph{#1}}

\begin{document}


\title{Poisson-Kingman Partitions\thanks{Research supported in part by N.S.F. Grants MCS-9404345 and DMS-0071448}}

\author{Jim Pitman
\\
\\
Technical Report No. 625 
\\
\\
Department of Statistics\\
University of California\\
367 Evans Hall \# 3860\\ 
Berkeley, CA 94720-3860
}

\date{October 23, 2002}

\def\shorttitle{Poisson-Kingman partitions}

\def\shortauthor{J. Pitman}

\maketitle

\thispagestyle{empty}


\abstract{
This paper presents some general formulas for random partitions 
of a finite set derived by Kingman's model of 
random sampling from an interval partition generated by subintervals whose
lengths are the points of a Poisson point process. These lengths can be also interpreted
as the jumps of a subordinator, that is an increasing process with stationary independent increments.
Examples include the two-parameter family of Poisson-Dirichlet models derived from the Poisson process of 
jumps of a stable subordinator.
Applications are made to the random partition generated by the lengths of excursions of a Brownian motion
or Brownian bridge conditioned on its local time at zero.
}

\keywords{exchangeable; stable; subordinator; Poisson-Dirichlet; distribution}

\newcommand{\Dir}{Dirichlet}
\newcommand{\fGamma}{\lambda}
\newcommand{\fSigma}{\lambda_\Sigma}
\newcommand{\giv}{\,|\,}
\newcommand{\Pdi}{P_{(i)}}
\newcommand{\str}{\tilde{\mu}}
\newcommand{\tilmu}{\tilde{\mu}}
\newcommand{\nqu}{\nu}
\newcommand{\pmk}{p_{[\kappa,\theta]}}
\newcommand{\DIR}{{\sc dirichlet}}
\newcommand{\dsal}{{\sc discrete stable}$(\alpha)$}
\newcommand{\PPP}{{\sc ppp}}
\newcommand{\GEM}{{\sc gem}}
\newcommand{\SPEC}{{\sc species}}
\newcommand{\PD}{{\sc pd}}
\newcommand{\PDth}{{\sc pd}$\,(\theta)$}
\newcommand{\DIRth}{{\sc dir}$\,(\theta Q )$}
\newcommand{\GEMth}{{\sc gem}$\,(\theta)$}
\renewcommand{\SS}{{\cal S}}
\newcommand{\DD}{{\cal D}}
\renewcommand{\AA}{{\cal A}}
\newcommand{\CC}{\mbox{${\mbox{\boldmath$M$}}$}}
\newcommand{\MM}{\mbox{${\mbox{\boldmath$M$}}$}}
\newcommand{\MC}{{M}}
\newcommand{\nj}{{jn}}
\renewcommand{\ni}{{in}}
\newcommand{\none}{{1n}}
\newcommand{\ntwo}{{2n}}
\newcommand{\FF}{{\cal F}}
\newcommand{\nnu}{{\nu}}
\newcommand{\RDS} {{\sc rd}$(S)$}
\newcommand{\RDN} {{\sc rd}$(\IN)$}
\newcommand{\RAM} {{\sc ram}}
\newcommand{\GRAM} {{\sc gram}}
\newcommand {\PP} {\mbox{\msbm P}}
\newcommand {\dagg} {\mbox{\dag}}
\renewcommand {\EE} {\mbox{\msbm E}}
\newcommand {\PTIL} {\mbox{$\tilde{\mbox{\boldmath$P$}}$}}
\newcommand {\TILP} {\mbox{$\tilde{\mbox{\boldmath$P$}}$}}
\newcommand {\BP} {\mbox{${\mbox{\boldmath$P$}}$}}
\newcommand {\bfone} {\mbox{${\mbox{\boldmath$1$}}$}}
\newcommand {\mm} {\mbox{${\mbox{\boldmath$m$}}$}}
\newcommand {\nnn} {\mbox{${\mbox{\boldmath$n$}}$}}
\newcommand {\nnnn} {(n_1, \ldots, n_k)}
\newcommand {\ppp} {\mbox{$p({\mbox{\boldmath$n$}})$}}
\newcommand {\Nns} {\mbox{${\mbox{\boldmath$N_n$}}$}}
\newcommand{\hf} {{\mbox{${\textstyle\frac{1}{2}}$}}}
\newcommand {\alth}{ {(\alpha, \theta)}}
\newcommand {\oth}{ {(0, \theta)}}
\newcommand {\althnu}{ {\alpha, \theta, \nu}}
\newcommand {\palth}{ {p_{\alpha, \theta} }}
\newcommand {\Palth}{ {\PP_{\alpha, \theta} }}
\newcommand {\Poth}{ {\PP_{0, \theta} }}
\newcommand {\malth}{ {(- \alpha, \theta)}}
\newcommand {\Pal}{ \mbox{$\cal{P}_{\alpha} $}}
\newcommand {\Qext}{ \mbox{$\cal{Q}_{\rm ext} $}}
\newcommand {\Q}{ \mbox{$\cal{Q} $}}
\newcommand{\PPp}{\mbox {${\cal P}$}}
\newcommand {\PPt}{ \mbox{$ {\cal P}_2 $}}
\newcommand {\PPb}{ \mbox{$\overline{ \cal P }_2 $}}
\newcommand {\calth}{ \mbox{$c_{\alpha, \theta}$}}
\newcommand {\la}{ {\lambda}}
\newcommand {\lai}{ {\lambda_i}}
\newcommand {\Mjn}{ \mbox{$M_{j,n}$}}
\newcommand {\Mjnt}{ \mbox{$M_{j,n_t}$}}
\newcommand {\Knt}{ \mbox{$K_{n_t}$}}
\newcommand {\LL}{ \mbox{$L$}}
\newcommand {\jt} {{j,t}}
\newcommand {\pjt}{ \mbox{$p_{j,t}$}}
\newcommand {\mujt}{ {\mu_{j,t} }}
\newcommand {\Mnjt}{ {M_{j,n_t} }}
\newcommand {\Mnj}{ {M_{j,n} }}
\newcommand {\Mnone}{ {M_{1,n} }}
\newcommand {\Mntwo}{ {M_{2,n} }}
\newcommand {\pjal}{ {p_{\alpha, j} }}
\newcommand {\palj}{ {p_{\alpha, j} }}
\newcommand {\pali}{ {p_{\alpha, i} }}
\newcommand {\cnt}{ \mbox{$c_{n,t}$}}
\newcommand {\Pla}{ \mbox{$P_{\la}$}}
\newcommand {\aaa} {{\kappa}}
\newcommand {\bbb} {{b}}
\newcommand {\ccc} {{c}}
\renewcommand {\ggg} {{\gamma}}
\newcommand {\alb} {{\alpha, \bbb}}
\newcommand {\LLth} {{L_{\theta}}}
\newcommand {\Qth} {{Q_{\theta}}}
\newcommand {\hth} {{h_{\theta}}}
\newcommand {\rhoth} {{g_{\theta}}}
\newcommand {\tot}{ \mbox{$\lambda_{\Sigma}$}}
\newcommand {\ttt} {{\tau}}
\newcommand {\Ptt} {{P_{\tau}}}
\newcommand {\rg} {{(\rho, g)}}
\newcommand {\Pg} {{P_g}}
\newcommand {\abc} {{(\alpha,b,c)}}
\newcommand {\Pabc} {{P_{(\alpha,b,c)}}}
\newcommand {\fal} {{f_{\alpha}}}
\newcommand {\Palo} {{\PP_{\alpha, 0}}}

\newcommand{\TT}{{\cal T}}
\newcommand{\GG}{{\cal G}}
\newcommand{\INk}{{[k]}}
\newcommand{\edge}{{-\!\!-}}
\newcommand{\iid}{{independent and identically distributed}}
\newcommand{\summ}[2]{{\Sigma_{#1}{#2}}}
\newcommand{\pnKN}{p(n,K,N)}
\newcommand{\Fnij}{F_{nN}^{ij}}
\newcommand{\lamnN}{\lambda_{nN}}
\newcommand{\lah}{\hat{\lambda}}
\newcommand{\vnN}[1]{V_{nN}^{(#1)}}
\newcommand{\vnf}{V_{nN}^{(4)}}
\newcommand{\nshv}{Y_{nN}}
\newcommand{\asn}{{\mbox{ as } n \te \infty }}
\newcommand{\nv}{{\#_{verts} }}
\newcommand{\TTN}{ {\cal T} _N }
\newcommand{\TTn}{ {\cal T} _n}
\newcommand{\TTk}{ {\cal T} _k }
\newcommand{\numt}[2]{\#_{#1} (#2)  }
\newcommand{\dkn}{\widehat{\#}_{n} (k)  }
\newcommand{\mupl}{\mu_+}
\newcommand{\errlap}[1]{\mu_+(#1)}
\newcommand{\apeq}{\simeq}
\newcommand{\novN}{\mbox{${n \over N }$} }
\newcommand{\kovN}{\mbox{${k \over N }$} }
\newcommand{\sspl}{v_+}
\newcommand{\rank}{\rho}
\newcommand{\deljmx}{W(m,j,x)}
\newcommand{\delj}[3]{W(#1,#2,#3)}
\newcommand{\numtp}[2]{\#_{#1} ^\perp (#2)  }
\newcommand{\npshapes}{\#_{shapes} ^{proper}  }
\renewcommand{\ss}{{s}}
\newcommand{\RT}{{\cal R}}
\newcommand{\RR}{{\cal R}}
\renewcommand{\NN}{{\bf N}}
\newcommand{\out}{\mbox{{\rm out}}}
\newcommand{\hF}{D}
\newcommand{\xell}{x}
\newcommand{\PINK}{\Pi_{KN}}
\newcommand{\PINJ}{\Pi_{JN}}
\newcommand{\PIN}{\Pi_{N}}
\newcommand{\PINN}{\Pi_{n,N}}
\newcommand{\KNP}{K_{N,p}}
\newcommand{\psnk}{p_{s,N}(n_1, \cdots, n_k)}
\newcommand{\ps}{p_s(n_1, \cdots, n_k)}
\newcommand{\Pisn}{\Pi_{s,n}}
\newcommand{\hPi}{\hat{\Pi}}
\newcommand{\tilP}{\tilde{P} }
\newcommand{\PPi}{{\cal P} _\infty }
\newcommand{\PPbar}{{\cal P} }
\newcommand{\Pn}{{\cal P} _n}
\newcommand{\Fn}{{\cal F} _n}
\newcommand{\Pnn}{{\cal P} _{[n]}}
\newcommand{\pisn}{\Pi_n(s)}
\newcommand{\simp}{\Delta}
\newcommand{\simb}{\bar{\Delta}}
\newcommand{\simd}{\nabla}
\newcommand{\lamda}{\lambda}
\newcommand{\ga}{\beta}
\newcommand{\cc}{C}
\newcommand{\simdb}{\bar{\nabla}}
\newcommand{\sidbar}{\bar{\nabla}}
\newcommand{\xx}{{\bf x }}
\renewcommand{\mm}{{\bf m }}
\newcommand{\nn}{{\bf n }}
\newcommand{\XX}{{\bf X }}
\newcommand{\YY}{{\bf Y }}
\newcommand{\VV}{{\bf V }}
\newcommand{\Qell}{{Q(\,\cdot\,|\,\ell) }}
\newcommand{\Qzero}{{Q(\,\cdot\,|\,0) }}
\newcommand{\sig}{{m}}
\newcommand{\sip}[2]{{  \simp_{#1}^{#2}  }} 
\newcommand{\sid}[2]{{  \simd_{#1}^{#2}  }} 
\newcommand{\tail}[2]{{  \mbox{{\sc tail}}( {#1}, {#2} )  }} 
\newcommand{\tai}[1]{{  \mbox{{\sc tail}}( {#1} ) }} 
\newcommand{\Psh}{\Pi_N^{\dagger} }
\renewcommand{\PD}{\mbox{{\sc pd}}}
\newcommand{\SBPD}{\mbox{{\sc sbpd}}}
\newcommand{\EPF}{\mbox{{\sc eppf}}}
\newcommand{\PGW}{\mbox{{\sc pgw}}}
\newcommand{\sbpd}{\mbox{{\sc sbpd}}}
\newcommand{\PDg}[1] {{\sc PD }( \hf | \hf {#1} )}
\newcommand{\unif}{uniform($N$)}
\newcommand{\tilXN}{{\tilde{\XX}^N (u)}}
\newcommand{\tXX}{{\tilde{\XX}}}
\newcommand{\sbfst}{\tilde{\Delta}_*^1}
\newcommand{\tilQl}{\tilde{Q}(\ell) }
\newcommand{\tilQ}{\tilde{Q}(\ell) }
\newcommand{\tilX}{\tilde{X}}
\newcommand{\Jelta}{\tilde{J}}
\newcommand{\VP}{\tilde{P}}
\newcommand{\RP}{\tilde{R}}
\newcommand{\VR}{\tilde{R}}
\newcommand{\uu}{r}
\newcommand{\freqs} {{\mbox{{\sc freqs}}}}
\newcommand{\ranked} {{\mbox{{\sc ranked}}}}
\newcommand{\tf}[1] {{T_{FIRST}^{(#1) } }}
\newcommand{\tl}[1] {{T_{LAST}^{(#1) } }}
\newcommand{\barPhi}{{\bar { \Phi } }}
\newcommand{\thf}{\mbox{$3 \over 2$ }}
\newcommand{\fth}{\mbox{$4 \over 3$ }}
\newcommand{\tilf}{\tilde{f}}
\newcommand{\scrhf}{{\scriptstyle{1 \over 2 }}}
\newcommand{\scrthf}{{\scriptstyle{3 \over 2 }}}
\newcommand{\scrtth}{{\scriptstyle{2 \over 3 }}}
\newcommand{\shf}{{\scriptstyle{1 \over 2 }}}
\renewcommand{\th}{\mbox{$1 \over 3 $}}
\newcommand{\qu}{\mbox{$1 \over 4 $}}
\newcommand{\PIBM}[1]{{\Pi_{#1} ^{BM}}}
\newcommand{\PIB}[1]{{\Pi_{#1}^{BM}}}
\newcommand{\PiB}{{\Pi^{B}}}
\newcommand{\PiBM}{{\Pi^{BM}}}
\newcommand{\PiBB}{{\Pi^{BB}}}
\newcommand{\rex}[1]{{\mbox{{\sc pd}}( \hf \giv #1)}}
\newcommand{\sbex}[1]{{\mbox{{\sc sbpd}}( \hf \giv #1)}}
\newcommand{\PTRN}[1]{\Pi_{#1}^{T_N}}
\newcommand{\PCRT}[1]{\Pi_{#1}^{T_\infty} }
\newcommand{\LMIN}[1]{L_{\min ( #1 ) } }
\newcommand{\LMAX}[1]{L_{\max ( #1 ) } }
\newcommand{\fhf}{{f_ {\mbox{${1 \over 2}$} }}}
\newcommand{\tilnu}{\tilde{\nu}}
\newcommand{\barp}{\bar{p}}
\newcommand{\QP}{Q}
\newcommand{\fm}{\Psi}
\newcommand{\mf}{\Psi}
\newcommand{\nnnt}{(n_1, \cdots, n_k \giv t ) }
\newcommand{\thoval}{\mbox{$\theta \over \alpha $}}
\newcommand{\beoval}{\mbox{$\beta \over \alpha $}}
\newcommand{\alo}{\alpha}
\newcommand{\tilV}{\tilde{V}}
\newcommand{\tilR}{{R}}
\newcommand{\gG}{{\gamma}}
\newcommand{\alal}{{\alpha,\alpha}}
\newcommand{\palthnnn}{{p_\alth \nnnn }}
\newcommand{\paltnnn}{{p_\alpha (n_1, \cdots, n_k  \giv t ) }}
\newcommand{\PKK}{{Poisson-Kingman}}
\newcommand{\PK}{\mbox{{\sc pk}}}
\newcommand{\sizeb}{\mbox{{\sc size-biased}}}
\newcommand{\sthf}{{stable$(\hf)$}}
\newcommand{\PDalth}{\mbox{{\sc pd$(\alpha ,\theta )$}}}
\newcommand{\PDalo}{\mbox{{\sc pd$(\alpha , 0 )$ }}}
\renewcommand{\PDth}{\mbox{{\sc pd$(\theta)$}}}
\newcommand{\givth}{{ \,|| \theta }}
\newcommand{\Calth}{{ C(\alpha, \theta }}
\renewcommand{\alth}{{ \alpha, \theta }}
\newcommand{\rhob}{{ \rho^{(b)}}}
\newcommand{\Pb}{{ P^{(b)}}}
\newcommand{\PO}{{ P_\alpha^0 }}
\newcommand{\EO}{{ E_\alpha^0 }}
\newcommand{\RRR}{{ \cal R }}
\newcommand{\fb}{{ f^{(b)}}}
\newcommand{\ells}{{ \ell }}
\newcommand{\gamb}{{ \gG^{(b)}}}
\newcommand{\aldiv}{{ \mbox{$\alpha$-{\sc diversity}}}}
\newcommand{\siml}{\mbox{$ \ \stackrel{\ell}{\sim}$ }}
\newcommand{\simn}{\mbox{$ \ \stackrel{n}{\sim}$ }}
\renewcommand{\sp}{\mbox{\sc l-span}}
\renewcommand{\span}{\mbox{\sc l-span}}
\newcommand{\vsp}{\mbox{\sc v-span}}
\newcommand{\dsp}{\mbox{\sc span}}
\newcommand{\dspst}{\mbox{\sc red-span}}
\newcommand{\tauN}{\tau_N }
\newcommand{\XnN}{X_{nN} }
\renewcommand{\ttt}{\hat{\tau}_{nN} }
\newcommand{\ST}[1]{ { \cal S  } _{#1} }
\renewcommand{\SS}[2]{ { \cal S  } _{#1, #2  } }
\newcommand{\siginv}{\sigma^{-1}}
\newcommand{\sigtil}{\tilde{\sigma}}
\newcommand{\prshape}{t^*}
\newcommand{\taunN}{\tau_{nN}}
\newcommand{\taumN}{\tau_{mN}}
\newcommand{\taun}{\tau_{n}}
\newcommand{\tauk}{\tau_{k}}
\newcommand{\pnn}[1]{ p_{nN}(#1) }
\newcommand{\ssigma}{\sigma}
\newcommand{\ttaun}{\tilde{\tau}_{n}}
\newcommand{\ttaunp}{\tilde{\tau}_{n+1}}
\newcommand{\ttaum}{\tilde{\tau}_{m}}
\newcommand{\ttaunN}{\tilde{\tau}_{nN}}
\newcommand{\epsnN}{\eps_{nN}}
\newcommand{\epst}{\eps_{3}}
\newcommand{\epsNN}{\eps_{NN}}
\newcommand{\nunN}{\nu_{nN}}
\newcommand{\munN}{\mu_{nN}}
\newcommand{\integ}[1]{\lceil #1 \rceil }
\newcommand{\delt}[1]{\delta_{#1}}
\newcommand{\nun}{\nu_{n}}
\newcommand{\etaN}{\eta_N}
\newcommand{\shape}{\tau^*}
\newcommand{\nunp}{\nu_{n+1}}
\newcommand{\rhonN}{\rho_{nN}}
\newcommand{\trho}{\tilde{\rho}}
\newcommand{\trhonN}{\tilde{\rho}_{nN}}
\newcommand{\hrhonN}{\hat{\rho}_{nN}}
\newcommand{\nutN}{\nu_{2N}}
\newcommand{\shnN}{{\tau}_{nN}^*}
\newcommand{\xinN}{\xi_{nN}}
\newcommand{\jnN}{J_{nN}}
\newcommand{\dd}{d}
\newcommand{\nedg}{\#_{edges}}
\newcommand{\nverts}{\#_{verts}}
\newcommand{\nends}{\#_{ends}}
\newcommand{\nend}{\#_{ends}}
\newcommand{\njun}{\#_{juncts}}
\newcommand{\kk}{k}
\newcommand{\stirl}[2]{ \mbox{ $ \left\{ \begin{array}{c} #1 \\ #2 \end{array}   \right\} $ }} 
\newcommand{\shapes}[1]{ { \cal S  } _{#1}^{*} }
\newcommand{\nsh}[2]{\# _{#1}^{*}(#2)}
\newcommand{\nshnmd}{\# _n^{*}(d, m )}
\newcommand{\nshnmn}{\# _n ^{*}(n ,m )}
\newcommand{\nshmtwo}{\# _n^{*}(2 ,m )}
\newcommand{\barc}[2]{\bar{C}_{#1,#2}}
\newcommand{\hnu}[2] { \hat{\nu}_{#1 #2} }

\newcommand{\RF}{\tilde{F}}
\newcommand{\cspdot}{\,\cdot\,}
\renewcommand{\Pb}{\PR^{(b)}}
\renewcommand{\nnn} {(n_1, \cdots, n_k  )}
\renewcommand{\PO}{\PR_\alpha^0}
\renewcommand{\EO}{\ER_\alpha^0}
\newcommand{\quart}{ { 1 \over 4 } }
\renewcommand{\PIB}[1]{{\Pi_{#1}}}
\newcommand{\bgiv}{ \,||\,}
\renewcommand{\rex}[1]{{\PD ( \hf \bgiv #1)}}
\newcommand{\RW}{ W }
\renewcommand{\ll}{ \lambda }
\newcommand{\xgamma}{ \xi }
\newcommand{\sdec}{{\cal P}^{\downarrow}}
\long\def\comment#1{{}\hfil\break}

\section{Introduction}
This paper
presents some 
general formulas for random partitions 
of a finite set derived by Kingman's model of 
random sampling from an interval partition generated by subintervals whose
lengths are the points of a Poisson point process.
Instances and variants of this model have found applications in the diverse 
fields of population genetics \cite{et95,grotespeed02}, 
combinatorics \cite{abt01,csp}, Bayesian statistics \cite{james02pp}, 
ecology \cite{mc65,en78}, statistical physics 
\cite{derrida81,derrida94,derrida97,ruelle87,tala01}, and
computer science \cite{janson01h}.

Section \ref{secstruc} recalls some general results for partitions obtained by sampling 
from a random discrete distribution.
These results are then applied in Section \ref{secpk} to 
the \PKK\ model. Section \ref{secops} discusses three basic operations on
\PKK\ models: scaling, exponential tilting, and deletion of classes.
Section \ref{secex} then develops formulas for specific examples
of \PKK\ models. Section \ref{twoparam} recalls the two-parameter family of Poisson-Dirichlet
models derived in \cite{py95pd2} from the Poisson process of jumps of 
a stable$(\alpha)$ subordinator for $0 < \alpha < 1$.
Section \ref{excns} reviews some results of \cite{ppy92,jp.bmpart,py92,py95pd2} relating the two-parameter family to the lengths of excursions of a Markov 
process whose zero set is the range of a stable subordinator of index $\alpha$.
Section \ref{sechf} provides further detail in
the case $\alpha = \hf$  which corresponds to 
partitioning a time interval by the lengths of excursions of a Brownian motion.
As shown in \cite{ap92,jpda97sac}, it is this \sthf\ model which 
governs the asymptotic distribution of partitions derived in various
ways from random forests, random mappings, and the additive coalescent.
See also \cite{bertoin00f,chassjan01} for further developments in terms of 
Brownian paths, and \cite{chasslou99,janson01h} for applications to hashing 
and parking algorithms.  This paper is a revision of the earlier 
preprint \cite{pitman95pk}. See \cite{csp} for a broader context and
further developments.

\section{Preliminaries}
\label{secstruc}
This section recalls some basic ideas from Kingman's theory of
exchangeable random partitions \cite{ki78b,ki82co}, as further developed in \cite{jp.epe}.
See \cite{jp96bl,csp} for more extensive reviews of these ideas and their
applications.
Except where otherwise specified, all random variables are assumed to
be defined on some background probability space $(\Omega, \FF, \PR)$,
and $\ER$ denotes expectation with respect to $\PR$.
Let $\IN:= \{1,2, \ldots \}$, let $F$ denote a random probability 
distribution on the line, and let $\Pi$ be a random partition of 
$\IN$ {\em generated by sampling from $F$}.
That is to say, two positive integers $i$ and $j$ are in the same block of 
$\Pi$ iff $X_i = X_j$, where conditionally given $F$ 
the  $X_i$ are independent and identically distributed according to $F$.
Formally, $\Pi$ is identified with the sequence $(\Pi_n)$, where $\Pi_n$ is the restriction of $\Pi$ to the finite set 
$\IN_n:= \{1, \ldots, n\}$. 
The distribution of $\Pi_n$ is such that for each particular partition $\{A_1, \cdots, A_k \}$ of $\IN_n$ with $\#(A_i ) = n_i$ for $ 1 \le i \le k $,
where $n_i \ge 1$ and $\sum_{i = 1}^k n_i = n$,
\eq
\lb{eppfdef}
\PR ( \Pi_n = \{A_1, \cdots, A_k \} ) = p(n_1, \cdots, n_k)
\en
for some symmetric function $p$ of sequences of positive integers, called the
{\em exchangeable partition probability function (EPPF)} of $\Pi$.
Conversely, Kingman \cite{ki78b,ki82co} showed that if $\Pi$ is an exchangeable
random partition of $\IN$,
meaning that the distribution of its restrictions $\Pi_n$ is of the 
form \re{eppfdef} for every $n$, for some symmetric function $p$, then $\Pi$ 
has the same distribution as if generated by sampling from some random probability distribution $F$.
Let $P_i$ denote the size of the $i$th largest atom of $F$. If 
$F$ is a random discrete distribution, then $\sum_i P_i = 1$ almost surely,
and $\Pi$ is said to have {\em proper frequencies} $(P_i)$.
In that case, let $\tilde{P}_j$ denote the size of the $j$th atom discovered in the process of random sampling. Put another way,
$\tilde{P}_j$ is the asymptotic frequency of the $j$th class of $\Pi$ when the classes are put in order of their 
least elements.
It is assumed now for simplicity that $P_i >0$ for all $i$ almost surely,
and hence $\tilde{P}_j > 0$ for all $j$ almost surely.
The sequence $(\tilde{P}_j )$ is a {\em size-biased permutation } of $(P_i)$.  That is to say, $\tilde{P}_j = 
P_{\pi_j}$ where for all finite sequences $(i_j, 1 \le j \le k)$ 
of distinct positive integers, the conditional probability 
of the event ($\pi_j = i_j$ for all $1 \le j \le k$) given $(P_1, P_2, 
\ldots)$ is
\eq
\lb{sbform}
P_{i_1} \, { P_{i_2}  \over 1 - P_{i_1}} \,
\cdots \, { P_{i_k}  \over 1 - P_{i_1} - \ldots - P_{i_{k-1}} } .
\en
The distribution of $\Pi_n$ is determined by the distribution of the sequence of ranked
frequencies $(P_i)$ through the distribution of the size-biased permutation 
$(\tilde{P}_j)$. To be precise, the EPPF $p$ in \re{eppfdef} is given by the formula \cite{jp.epe}
\eq
\lb{genepf}
p(n_1, \cdots, n_k  ) = 
\ER  \left[ \prod_{i=1}^k \tilP_i^{n_i - 1} \prod_{i=1}^{k-1} \left(1 - \sum_{j=1}^i \tilP_j \right)\right] .
\en
Alternatively \cite{jp96bl}
\eq
\lb{permsum}
p \nnn = \sum_{(j_1, \ldots, j_k )} \ER \prod_{ i = 1}^k P_{j_i}^{n_i}
\en
where $(j_1, \ldots, j_k )$ ranges over all permutations of $k$ positive 
integers, and the same formula holds with $P_{j_i}$ replaced by $\tilP_{j_i}$.
For each $n = 1,2, \cdots $ the EPPF $p$,
when restricted to $(n_1, \cdots, n_k )$ with $\sum_i n_i = n$,
determines the distribution of $\Pi_n$.  
Since $\Pi_{n}$ is the restriction of $\Pi_{n+1}$ to $\IN_n$,
the EPPF is subject to the following sequence of {\em addition rules}
\cite{jp.epe}: for $k = 1,2, \ldots$
\eq
\lb{addr}
p \nnn = \sum_{j = 1}^k p( \ldots, n_j+1, \ldots) + 
p (n_1, \ldots, n_k,1) 
\en
where $( \ldots, n_j+1, \ldots)$ is  derived from
$(n_1, \ldots, n_k)$ by substituting $n_j+1$ for $n_j$.
The first few rules are
\eq
\lb{cons1}
1 = p(1) = p(2) + p(1,1)
\en
\eq
\lb{cons2}
p(2) = p(3) + p(2,1)  ; ~~~~~p(1,1) = 2 p(2,1) + p(1,1,1)
\en
where $p(2,1) = p(1,2)$ by symmetry of $p$.
Let $\mu ( q )$ denote the $q$th moment of $\tilP_1$:
\eq
\lb{cons4}
\mu ( q ) := \ER [ \, \tilP_1 ^{q} \, ]  = \int_0^1 p^q \, \tilnu(dp).
\en
where $\tilnu$ denotes the 
distribution of $\tilP_1$ on $(0,1]$.
Following Engen \cite{en78}, call $\tilnu$
the {\em structural distribution} associated with an
random discrete distribution whose size-biased permutation is $(\tilP_j)$,
or with an exchangeable random partition $\Pi$ whose sequence of 
class frequencies is $(\tilP_j)$.
The special case of \re{genepf} for $k=1$ and $n_1 = n$
is
\eq
\lb{cons3}
p(n) = \ER [ \, \tilP_1 ^{n-1} \, ] = \mu ( n - 1)
~~~~~(n = 1, 2, \cdots ).
\en
From \re{cons1}, \re{cons2}, and \re{cons3} the following 
values of the EPPF are also determined by the first two moments of the
structural distribution:
\eq
\lb{cons5}
p(1,1) = 1 - \mu(1) ;~~  p(2,1) = \mu(1) - \mu(2) ;~~ p(1,1,1) = 1 - 3 \mu(1) + 2 \mu(2).
\en
So the distribution of the random partition of $\{1,2,3\}$ induced by $\Pi$ 
with class frequencies $(\tilP_i)$ is determined by the first
two moments of the structural distribution of $\tilP_1$.
It is not true in general that the EPPF
is determined for all $(n_1, \cdots, n_k)$ by the
structural distribution, because it is possible to construct different
distributions for a sequence of ranked frequencies which have
the same structural distribution.

Continuing to suppose that $(P_i)$ is the sequence of ranked atoms of
a random discrete probability distribution,
and that $(\tilP_j )$ is a size-biased permutation of $(P_i)$,
for an arbitrary non-negative measurable function $f$, there is the 
well known formula 
\eq
\lb{eform}
\ER \left[ \sum_i f ( P_i ) \right] = 
\ER \left[ \sum_j f ( \tilP_j ) \right] = \ER \left[ {f( \tilde{P}_1) 
\over \tilP_1 } \right]
= \int_0^1 {f (p ) \over p } \tilnu ( dp ).
\en
This formula shows that
the structural distribution $\tilnu$
encodes much information about the entire sequence
of random frequencies.
Taking $f$ in \re{eform}
to be the indicator of a subset $B$ of $(0,1]$, the quantity in
\re{eform} is $\nu(B) = \int_B p^{-1} \tilnu (dp)$. This measure $\nu$
is the mean intensity measure of the point process with a point at each 
$P_i\in (0,1]$. 
For $x > \hf$ there can be at most one $P_i > x$, so the
structural distribution $\tilnu$ determines the distribution of 
$P_1 = \max_j \tilP_j $ on $(\hf ,1]$ via the formula
\eq
\PR ( P_1 > x ) = \nu(x,1] = \int_{(x,1] }p ^{-1} \tilnu(dp) ~~~~~~(x > \hf) .
\en
Typically, formulas for $\PR ( P_1 > x )$ get progressively more complicated
on the intervals $(\th, \hf]$, $(\qu, \th], \cdots$.
See for instance \cite{per91,py95pd2}.

A random variable of interest in many applications is the sum of $m$th
powers of frequencies
$$
S_m:= \sum_{i= 1}^\infty P_{i}^{m} = \sum_{j= 1}^\infty \tilP_{j}^{m}
~~~~~~~~~(m = 1,2, \ldots)
$$
where it is still assumed that $S_1 = 1$ almost surely.
Let $\pi := \{A_1, \cdots, A_k \}$ be some particular partition of
$\IN_n$ with $\#(A_i ) = n_i$ for $ 1 \le i \le k $, and
consider the event $(\Pi_n \ge \pi)$, meaning that each block of 
$\Pi_n$ is some union of blocks of $\pi$. Then it is easily shown that
\eq
\lb{qn2}
\PR(\Pi_n \ge \pi ) = \ER \left[ \prod_{i=1}^k S_{n_i} \right] = \sum_{j = 1}^k ~~
\sum_{ \{ B_1, \ldots, B_j \} } p( n_{B_1}, \ldots, n_{B_j})
\en
where the second sum is over partitions $\{ B_1, \ldots, B_j \}$ 
of $\IN_k$, and $n_B:= \sum_{i \in B} n_i$.
In particular, for $n_i \equiv m$ this gives an expression for the
$k$th moment of $S_m$ for each $k = 1,2, \ldots$:
\eq
\lb{qn22}
\ER \left[ S_{m}^k \right] = 
\sum_{j = 1}^k ~{1 \over j!} \, \sum_{ (k_1, \ldots, k_j ) } { k ! \over k_1! \cdots k_j ! } ~p( m k_1, \ldots, m k_j ) 
\en
where the second sum is over all sequences of $j$ positive integers
$(k_1, \ldots, k_j )$ with $k_1+ \cdots + k_j = k$.
Thus the EPPF associated with a random discrete distribution 
directly determines the positive integer moments of the power sums $S_m$, 
hence the distribution of $S_m$, for each $m$.

\section{The Poisson-Kingman Model}
\label{secpk}

Following McCloskey \cite{mc65}, Kingman \cite{ki75},
Engen \cite{en78}, Perman-Pitman-Yor \cite{per91,ppy92,py95pd2},
consider the ranked random discrete distribution $(P_i) := (J_i/T)$
derived from an inhomogeneous Poisson point process of random lengths 
$J_1 \ge J_2 \ge \cdots \ge 0$ by normalizing these lengths by
their sum $T := \sum_{i=1}^\infty J_i$.
So it is assumed that the number $N_I$ of $J_i$ that fall in an interval $I$ is 
a Poisson variable with mean $\Lambda(I)$, for some \Lev\ measure $\Lambda$
on $(0,\infty)$,
and the counts $N_{I_1}, \cdots , N_{I_k}$
are independent for every finite collection of disjoint intervals
$I_1, \cdots , I_k$. 
It is also assumed that 
$$\int_0^1 x \Lambda(dx) < \infty \mbox{ and } \Lambda[1,\infty) < \infty$$
to ensure that $\PR(T < \infty ) = 1$.
The sequence $(P_i)$ may be regarded 
as a random element of the space $\sdec$ of decreasing sequences of positive real numbers with sum $1$.
Throughout this section, the following further assumption
is made to ensure that various conditional probabilities can
be defined without quibbling about null sets:

\noindent
{\bf Regularity assumption.}
{\em The \Lev\ measure $\Lambda$ has a density $\rho(x)$ such that the distribution of $T$ 
is absolutely continuous with density 
$$
f(t) := \PR (T\in dt)/dt
$$
which is strictly positive and continuous on $(0,\infty)$. }

Note that the regularity assumption implies the total mass of the
\Lev\ measure is infinite:
\eq
\lb{infy}
\int_0^\infty \rho(x) dx = \infty  .
\en
The results described below also have weaker forms for a
\Lev\ density $\rho(x)$ just subject to \re{infy}, with appropriate
caveats about almost everywhere defined conditional probabilities.

It is well known that $f$ is uniquely determined by $\rho$ 
via the Laplace transform 
\eq
\lb{lkf}
\ER ( e^{-\lambda T } ) = \int_0^\infty e^{-\lambda x } f(x) dx = \exp [ - \psi (\lambda )]
~~~~~~~~(\la \ge 0)
\en
where, 
according to the \LK\ formula,
\eq
\lb{psilam}
\psi(\lambda)= \int_0^\infty ( 1 - e^{- \lambda x } ) \rho(x) dx .
\en
Alternatively,
$f$ is the unique solution of the following integral equation,
which can be derived from \re{lkf} and \re{psilam} by differentiation
with respect to $\la$:
\eq
\lb{inteq}
f(t) = \int_0^t  \rho(v ) f(t - v ) { v \over t } dv   .
\en
Let $(\tilP_j)$ be a size-biased permutation of the normalized lengths $(P_i) := (J_i/T)$
and let $(\Jelta_j) = (T \tilP_j)$ be the
corresponding size-biased permutation of the ranked lengths $(J_i)$.
Then  \re{inteq} admits the following probabilistic interpretation \cite{mc65,ppy92}: 
\eq
\lb{jelta}
\PR ( \Jelta_1 \in dv, T \in dt ) = \rho(v ) dv f(t - v ) dt { v \over t } .
\en
This can be understood as follows.
The left side of \re{jelta} is the probability that among the Poisson
lengths there is some length in $dv$ near $v$, and
the sum of the rest of the lengths falls in an interval of length $dt$ near $t-v$, and finally that 
the interval of length about $v$ is the one picked by length-biased sampling.
Formally, \re{jelta} is justified by
the description of a Poisson process in terms of its Palm measures \cite{ppy92}.

The following two Lemmas are read from \cite[Theorem 2.1]{ppy92}.
The first Lemma is immediate from \re{jelta}, and the second
is obtained by a similar Palm calculation.

\begin{lmm}{ppy} {\em \cite{ppy92}}
For each $t > 0$ the formula
\eq
\lb{strucha}
\tilde{f}( p \giv t )
:= p t \, \rho ( p t )  \,
{f ( \barp t ) \over f ( t ) }   
~~~~~~(0 < p < 1 ; ~ \barp := 1 - p ),
\en
where $\rho$ is the density of the \Lev\ measure of $T$ and $f$ is the probability
density of $T$,
defines a function of $p$ which is a probability density on $(0,1)$.
This is the density of the structural distribution of
$\tilP_1 := \Jelta_1/T$ given $T = t$: 
\eq
\lb{stru2}
\PR ( \tilP_1 \in dp \giv t )  = \tilde{f}( p \giv t ) dp   ~~~~~~(0 < p < 1 ).
\en
\end{lmm}

\begin{lmm}{ppy2} {\em \cite{ppy92}}
For $j = 0,1,2, \cdots $ let 
\eq
\lb{tjs}
T_j := T - \sum_{k = 1}^j \Jelta_k = \sum_{k = j+1}^\infty \Jelta_k 
\en
which is the total length remaining after removal of the first
$j$ Poisson lengths $ \Jelta_1, \ldots, \Jelta_j$ chosen by length-biased sampling. 
Then the family of densities {\em \re{strucha}} on $(0,1)$,
parameterized by $t > 0$, provides the conditional density of the random variable
$$
G_{j+1} := {\Jelta_{j+1} \over T_j } = { \tilP_{j+1} \over \tilP_{j+1} + \tilP_{j+2} + \cdots }
$$
given $T_0, \cdots, T_j$ via the formula
\eq
\lb{lkf2}
\PR \left( G_{j+1} \in dp \left|  T_0, \cdots, T_j \right) \right.= \tilde{f}(p | T_j ) \,dp
~~~~~~~(0 < p < 1).
\en
\end{lmm}

Lemma \ref{ppy2} provides an explicit construction of a regular conditional
distribution for $(\tilP_j)$ given $T=t$ for arbitrary $t>0$. This
conditional distribution of $(\tilP_j)$ given $T=t$ determines corresponding
conditional distributions for the $\sdec$-valued 
ranked sequence $(P_i)$ and for an associated random
partition $\Pi$ of $\IN$.

\begin{dfn}{pkk}{\em
The distribution of $(P_i) := (J_i/T)$ on $\sdec$ determined
by the ranked points $J_i$ of a Poisson process with \Lev\ density $\rho$
will be called the {\em Poisson-Kingman distribution with \Lev\ density $\rho$},
and denoted $\PK(\rho)$.
Denote by \PK$(\rho \giv t)$ the regular conditional distribution of 
$(P_i)$ given $(T = t )$ constructed above. For a probability
distribution $\gG$ on $(0,\infty)$,
let 
\eq
\PK(\rho , \gG ) := \int_0^\infty \PK(\rho \giv t) \gG(dt)
\en
be the distribution on $\sdec$ obtained 
by mixing the \PK$(\rho \giv t)$ with respect to $\gG(dt)$.
Call \PK$(\rho , \gG )$ the {\em Poisson-Kingman distribution 
with \Lev\ density $\rho$ and mixing distribution $\gG$}.
}
\end{dfn}
Note that
\PK$(\rho \giv t)$ = \PK$(\rho , \delta_t )$, where $\delta_t$ is 
a unit mass at $t$, and that \PK$(\rho )$ = \PK$(\rho , \gG )$
for $\gG(dt) = f(t) dt$.
A formula for the joint density of $(P_1, \cdots, P_n)$ for
$(P_i)$ with \PK$(\rho \giv t)$ distribution was obtained by
Perman \cite{per91}
in terms of the joint density $p_1 (t,x)$
of $T$ and $J_1$. This function can be described in terms of $\rho$ and $f$ as the solution of an integral equation \cite{per91}, 
or as a series of repeated integrals \cite{py95pd2}.
But this formula will not be used here.

For a probability distribution $Q$ on $\sdec$,
such as $Q = \PK(\rho, \gG)$, a random partition $\Pi$ of $\IN$ will
be called a {\em $Q$-partition} if $\Pi$ is an exchangeable random partition
of $\IN$ whose ranked class frequencies are distributed according to $Q$.
Immediately from Definition \ref{pkk}, the structural distribution of a \PK$(\rho, \gG)$-partition $\Pi$ of $\IN$, 
that is the distribution on $(0,1)$ of the frequency $\tilP_1$ of the class of $\Pi$ containing $1$,
has density
\eq
\lb{structural}
\PR (\tilP_1 \in dp ) /dp = \int_0^\infty \tilde{f}(p \giv t ) \gG(dt) ~~~~~(0 < p < 1)
\en
where $\tilde{f}(p \giv t )$ given by \re{strucha} is the density of the structural 
distribution of $\tilP_1$ given $T=t$ in the basic Poisson construction.
Similarly, the EPPF of $\Pi$ is
\eq
\lb{poisepf}
p(n_1, \cdots, n_k) = \int_0^\infty p(n_1, \cdots, n_k \giv t) \gG(dt)
\en
where $p(n_1, \cdots, n_k \giv t)$, the EPPF of a \PK$(\rho \giv t )$-partition, is determined as follows:

\begin{thm}{prp33}
The EPPF of a \PK$(\rho \giv t )$-partition is given by the formula
\eq
\lb{sp.18a}
p( n_1, \cdots, n_k \giv t )
=  t ^ {k-1} \int _0 ^ {1} p ^{ n + k -2 } I( n_1, \cdots, n_k ; tp ) 
\tilde{f}(p \giv t ) dp
\en
where 
$n = \sum_1^k n_i$, $I( n; v) = 1$ if $k = 1$ and $n_1 = n$,
and for $k = 2, 3 , \ldots$
\eq
\lb{sp.18b}
I(n_1, \cdots , n_k ; v) := {1 \over \rho(v) } \int _{ {\cal S}_k} \left[ \prod_{i = 1}^k \rho ( v u_i ) u_i ^{n_i}\right]
du_1 \cdots d u _{k-1}
\en
where ${\cal S}_k$ is the simplex
$\{ (u_1, \ldots , u_k ): u_i \ge 0 \mbox{ and } u_1 + \cdots + u_k = 1 \}$.
\end{thm}
\proof
In view of the formula \re{strucha} for $\tilde{f}(p \giv t )$,
the formula \re{sp.18a} is obtained from formula \re{tmn} in
the following Lemma by dividing by $f(t) dt$, letting $p = \sum_i x_i / t$,
and integrating out with respect to $p$ and to $u_i = x_i / ( p t )$ for $1 \le i \le k-1$.
\endpf

A change of variables gives the following variant of formula \re{sp.18a}, 
whose connection to the next lemma is a bit more obvious:
\eq
\lb{sp.18aa}
p( n_1, \cdots, n_k \giv t )
=  \int _0 ^ {t} dv { f( t -v ) \over t^n f(t)  } v^{n+k - 1}
I (n_1, \ldots, n_k;v)  \rho(v).
\en
\begin{lmm}{lmpalm}
Let $\Pi_n$ be the restriction to $\IN_n$ of a $\PK(\rho)$ partition $\Pi$ whose
class frequencies (in order of least elements) are $\tilP_j = \Jelta_j/T$, where $T = \sum_{j} \Jelta_j$
has density $f$, and the lengths $\Jelta_j$ are the points of a Poisson process of lengths with intensity 
$\rho$, in length-biased random order.
Then for each partition $\{A_1, \cdots, A_k \}$ of $\IN_n$ such that
$\#(A_i ) = n_i $ for $ 1 \le i \le k $,
\eq
\lb{thisprob}
\PR ( \Pi_n = \{A_1, \cdots, A_k \} , \Jelta_i \in dx_i \mbox{ for } 1 \le i \le k, T \in dt )
\en
\eq
\lb{tmn}
= t^{-n} \, f ( t - \mbox{$\sum_{i=1}^k x_i$} )  \, dt \, \prod_{i=1}^k \rho ( x_i )  x_i^{n_i} \, d x _i   .
\en
\end{lmm}
\proof
This can be derived by evaluation of the expectation \re{genepf} for the joint distribution of
$\tilP_1, \ldots, \tilP_k$ given $T= t$ determined by Lemma \ref{ppy2}. 
Alternatively, there is 
the following more intuitive argument, which can be made rigorous using 
the characterization of Poisson process by its a Palm measures, as 
in \cite{py92,ppy92}. 
Let $\Pi$ be constructed as in \cite{jp.bmpart} using random 
intervals $I_i$ 
laid down on $[0,T]$
in some arbitrary random order, where the lengths $J_i:= |I_i|$ are the ranked
points of the Poisson process with intensity $\rho(x)$, and $T = \sum_i J_i$.
Let $U_1, U_2, \cdots $ be i.i.d. uniform on $(0,1)$ independent
of this construction.
Let $\Pi$ be the partition of $\IN$ generated
by the random equivalence relation $n \sim m $ iff either $n=m$
or $TU_n$ and $TU_m$ fall in the same interval $I_i$ for some $i$.
Then by construction, $\Pi$ is a $\PK(\rho)$ partition.
For the event in \re{thisprob} to occur,

(i) there must be some Poisson point in $dx_i$ for each $1 \le i \le k$, and

(ii) given (i), the sum of the rest of the Poisson points must fall in an interval of length $dt$ near $t - \sum_{i = 1}^k  x_i$, and

(iii) given (i) and (ii), for each $ 1 \le i \le k $ and each $ m \in A_i $
the sample point $T U_m$ must fall in the interval of length $x_i$.

The infinitesimal probability in \re{thisprob} therefore equals
\eq
\left( \prod_{i=1}^k \rho ( x_i )  \, d x _i  \right) \, f ( t - \mbox{$\sum_{i=1}^k x_i$} ) \, dt \,
\prod_{i =1}^k \left( {x_i \over t } \right) ^{n_i} 
\en
which rearranges as \re{tmn}.
\endpf

The formula \re{sp.18a} expresses
$p( n_1, \cdots, n_k \giv t )$ as the expectation of a function
of $\tilP_1$ given $T = t$, where the function depends on $t$ and
$n_1, \cdots, n_k$. 
Because some values of an EPPF can always
be expressed as moments of $\tilP_1$, as in \re{cons4} and \re{cons5},
it seems natural to try to express an EPPF similarly whenever possible.
This idea serves as a guide to simplifying calculations in a number of
particular cases treated later. 
The integrations in \re{sp.18a} and \re{sp.18b} are essentially convolutions, which can be expressed or evaluated in various ways. 
Consider for instance
the length $T_k:= T - \sum_{i=1}^k \Jelta_i$ which remains after removal of the first
$k$ lengths discovered by the sampling process. Then the formula of Lemma \ref{lmpalm} can be
recast as
\eq
\lb{thisprob2}
\PR ( \Pi_n = \{A_1, \cdots, A_k \} , \Jelta_i \in dx_i \mbox{ for } 1 \le i \le k, T_k \in dv )
\en
\eq
\lb{tmn1}
= (v + \mbox{$\sum_{i=1}^k x_i$} )  ^{-n} \, f(v) dv \, \prod_{i=1}^k \rho ( x_i )  x_i^{n_i} \, d x _i  
\en
which yields the following integrated forms of \re{sp.18a}: 

\begin{crl}{crlpk}
The EPPF of a \PK$(\rho)$-partition is given by the formula
\eq
\lb{eppf1}
p( n_1, \cdots, n_k )
= \int_0^\infty \cdots
\int_0^\infty \, { f(v) dv \, \prod_{i=1}^k \rho ( x_i )  x_i^{n_i} \, d x _i  \over
(v + \mbox{$\sum_{i=1}^k x_i$} )  ^n }
\en
where $n:= \sum_{i = 1}^k n_i$, or again by
\eq
\lb{eppf2}
p( n_1, \cdots, n_k )
= 
{ (-1)^{n-k} \over \Gamma (n ) } \int_0^\infty \, \la^{n-1} \, d \la e^{- \psi(\la) } \prod_{i = 1}^k \psi_{n_i} ( \la )
\en
where $\psi(\lambda):= \int_0^\infty ( 1 - e^{- \lambda x } ) \rho(x) dx$
is the Laplace exponent as in {\em \re{psilam}}, and 
\eq
\lb{psimla}
\psi_m(\la) := {d ^m \over d \la^m } \psi (\la) = (-1)^{m-1} \int_0^\infty x^{m} e^{-\la x }  \rho(x) dx ~~~
(m = 1,2, \ldots).
\en
\end{crl}
\proof
Formula \re{tmn1} yields \re{eppf1} by integration, and \re{eppf2} follows after applying the formula
$b^{-n} = \Gamma(n)^{-1} \int_0^\infty \la^{n-1} e^{- \la b} d \la$ to $b = v + \sum_{i=1}^k x_i$.
\endpf
These integrated forms \re{eppf1} and \re{eppf2} also hold more generally, with $f(v) dv$ replaced by $\PR (T \in dv)$,
and $\rho(x) dx$ replaced by the corresponding \Lev\ measure on $(0,\infty)$, assuming only that the \Lev\ measure has infinite 
total mass.

Provided $\ER (e^{\eps T}) < \infty$ for some $\eps >0$, the Laplace exponent $\psi$ can be expanded in a neighbourood
of $0$ as
$$
\psi( \la ) =  - \sum_{m = 1}^\infty  { \kappa_m \over m ! } ( - \la )^m
$$
where the {\em cumulants} $\kappa_m$ of $T$ are the moments of the \Lev\ measure
$$
\kappa_m = (-1)^{m-1} \psi_m(0) = \int_0^\infty x^{m} \rho(x) dx  .
$$
Then for each partition $\{A_1, \cdots, A_k \}$ of $\IN_n$ such that
$\#(A_i ) = n_i $ for $ 1 \le i \le k $,
Lemma \ref{lmpalm} yields the formula
\eq
\lb{xxpp}
\PR ( \Pi_n = \{A_1, \cdots, A_k \} , T \in dt ) = t^{-n} \PR ( T + \Sigma_{i = 1}^{k} J_{i,n_i} \in dt)
\prod_{ i = 1}^k \kappa_{n_i}
\en
where $J_{i,n_i}$ denotes a random length distributed according to the \Lev\ density tilted by
$x^{n_i}$:
$$
\PR ( J_{i,n_i} \in dx ) = \kappa_{n_i}^{-1} \rho(x) x^{n_i} \, dx
$$
and $T$ and the $J_{i,n_i}$ for $1 \le i \le k$ are assumed to be independent.
If $f_{n_1, \ldots, n_k}(t)$ denotes the probability density of $T + \Sigma_{i = 1}^{k} J_{i,n_i}$, then 
formula \re{sp.18a} for the EPPF of a \PK$(\rho \giv t )$-partition 
can be rewritten
\eq
\lb{sp.18ab}
p( n_1, \cdots, n_k \giv t ) = { f_{n_1, \ldots, n_k}(t)  \over t^{n} f(t) } 
\prod_{ i = 1}^k \kappa_{n_i}
\en
and formula \re{eppf1} for the EPPF of a \PK$(\rho)$-partition becomes
\eq
\lb{sp.19}
p( n_1, \cdots, n_k) = \ER \left[ ( T + \Sigma_{i = 1}^{k} J_{i,n_i} )^{-n} \right] \prod_{ i = 1}^k \kappa_{n_i} .
\en
See also James \cite{james02pp} for closely related formulas, with 
applications to Bayesian non-parametric inference.

\section{Operations}
\label{secops}

Later discussion of specific examples of Poisson-Kingman partitions will be guided 
by a number of basic operations on \Lev\ densities $\rho$ and their associated
families of partitions.

\subsection{Scaling}
By an obvious scaling argument, the \PK$(\rho)$ and
\PK$(\rho')$ distributions are identical whenever
$\rho'(x) = b \rho( b x )$ is a rescaling of $\rho$ for some $b > 0$.
The converse is less obvious, but true \cite[Lemma 7.5]{py92}.

\subsection{Exponential tilting}
It is elementary that if $\rho$ is a \Lev\ density,
corresponding to a density $f$ for $T$, and $b$ is a real number such that $\psi(b)$
defined by \re{psilam} is finite, then
\eq
\rhob (x ) = \rho (x ) e^ {- b x }
\en
is also a \Lev\ density, and the corresponding density of $T$
is
\eq
\lb{fb}
\fb ( t ) =  f( t ) \, e^{\psi (b ) - b t} 
\en
It is also well known \cite[Proposition 2.1.3]{ksor97} that if $\Pb$ denotes the
probability distribution governing the Poisson
set up with \Lev\ density $\rhob$ 
then \re{fb} extends to the absolute continuity relation
\eq
\lb{fb1}
{ d \Pb \over d \PR ^{(0)}}  = e^{\psi (b ) - b T}  .
\en
This relation is equivalent to a combination of \re{fb} and the following
identity, which can also be verified using the construction of Lemma \ref{ppy2}:
\eq
\PK ( \rhob \giv t ) = \PK ( \rho \giv t ) \mbox{ for all }  t > 0 .
\en
Consequently 
\eq
\PK ( \rhob , \gamma ) = \PK ( \rho, \gamma )
\en
for every $\gamma$. In particular, the distribution on $\sdec$
derived from the unconditioned Poisson model with \Lev\ density $\rhob$ is
\eq
\lb{pktilt}
\PK ( \rhob ) = \PK ( \rho, \gamb )
\en
where $\gamb$ is the $\Pb$ distribution of $T$, that is $\gamb(dt)  = \fb(t) dt$
for $\fb$ as in \re{fb}.
It can also be shown that if $\rho'$ and $\rho$ are two regular \Lev\ densities
such that $\PK ( \rho' ) = \PK ( \rho, \gG )$ for some $\gG$, then 
$\rho' = \rhob$ and $\gG = \gamb$ for some $b$.

\subsection{Deletion of Classes.}
The following proposition, which generalizes a result of \cite{ppy92},
provides motivation for study of \PK$(\rho, \gG)$-partitions
for other distributions $\gG$ besides $\gG(dt)=f(t)dt$ corresponding to
the unconditioned Poisson set up, and $\gG = \delta_t$ corresponding
to conditioning on $T= t$.
Given a random partition $\Pi$ of $\IN$ with infinitely many classes, for each 
$k = 0, 1, \cdots$ let $\Pi_k$ be the
partition of $\IN$ derived from $\Pi$ by 
{\em deletion of the first $k$ classes}, an operation made precise
as follows.
First let $\Pi_k'$ be the restriction of $\Pi$ to $H_k:= \IN - G_1 - \cdots - G_k$
where $G_1, \cdots G_k$ are the first $k$ classes of $\Pi$
in order of least elements, then derive $\Pi_k$ on $\IN$ from
$\Pi_k'$ on $H_k$ by renumbering the points of $H_k$ in increasing order.

\begin{prp}{prptk}
Let $\Pi$ be a \PK$(\rho, \gG)$-partition of $\IN$, 
and let $\Pi_k$ be derived from $\Pi$ be deletion of its first $k$ classes.
Then $\Pi_k$ is a \PK$(\rho, \gG_k)$-partition of $\IN$,
where $\gG_k = \gG Q^k$ for $Q$ the Markov transition operator on $(0,\infty)$
$$
Q(t, dv ) = \rho(t-v) (t-v) t^{-1} f(v) 1 ( 0 < v < t) dv  .
$$
In particular, if $\Pi$ is a \PK$(\rho)$ partition of $\IN$, then
$\Pi_k$ is \PK$(\rho, \gG_k)$-partition of $\IN$,
where $\gG_k$ is the distribution of $T_k$,
the total sum of Poisson lengths $T$ minus the sum of the first $k$ lengths discovered by a process of length-biased sampling,
as in {\em \re{tjs}}.
\end{prp}
\proof
According to a result of \cite{ppy92} which is implicit in Lemma \ref{ppy2},
the sequence $(T_k)$ is Markov with stationary transition probabilities
given by $Q$. The conclusion follows from this observation, the
construction of \PK$(\rho ; \gG )$,
and the general construction of an exchangeable partition of $\IN$ conditionally given its class frequencies \cite{jp.epe}.

\section{Examples}
\label{secex}

\subsection{The one-parameter Poisson-Dirichlet distribution.}
Following Kingman \cite{ki75}, for the particular choice
\eq
\lb{gamma}
\rho(x) = \theta \, {x}^{-1} \, e ^ {-b x }
\en
where $\theta > 0$ and $b > 0$, corresponding to $T$ with the
gamma$(\theta,b)$ density
\eq
\lb{gammadens}
f(t) = {b^\theta \over \Gamma(\theta) } t^{\theta -1} e^{ - b t },
\en
the \PK$(\rho)$ distribution is
the {\em Poisson-Dirichlet distribution with parameter $\theta$}, abbreviated \PDth.
Note the lack of dependence on the inverse scale parameter $b$.
The well known fact the structural distribution of \PDth\ is beta$(1,\theta)$ follows
immediately from \re{strucha}.
It follows easily from any one of the previous general formulas 
\re{sp.18a}, \re{eppf1}, \re{eppf2} or \re{sp.19}, 
that the EPPF of a \PDth-partition $\Pi = (\Pi_n)$ is given by the formula 
\eq
\lb{esf}
p_\theta (n_1, \cdots, n_k) =  {\theta^{k} \Gamma(\theta) \over \Gamma(\theta + n) }
\prod_{i = 1}^k (n_i - 1)! 
~~~~~~~(n = \sum_{i = 1}^k n_i ) .
\en
This is a known equivalent \cite{king82,jp.epe} of the Ewens sampling formula \cite{ew72,et95} for the joint 
distribution of the number of blocks of $\Pi_n$ of various sizes.
It is also known \cite{ppy92,py92} that the 
following conditions on $\rho$ are equivalent:

(i) $\rho$ is of the form \re{gamma}, for some $b>0, \theta >0$;

(ii) \PK$(\rho \giv t)= $\PK$(\rho)$ for all $t > 0$;

(iii) \PK$(\rho)= $\PDth\ for some $\theta > 0$.

(iv) a \PK$(\rho)$-partition has EPPF of the form \re{esf} for some $\theta > 0$.

\noindent
See also \cite{abt01,ki93} for further properties and applications of \PDth.
\subsection{Generalized gamma}
After the one-parameter Poisson-Dirichlet family, the next simplest
\Lev\ density $\rho$ to consider is
\eq
\lb{gamma1}
\rho_{\alpha,c,b}(x) = c  \, {x}^{- \alpha - 1} \, e ^ {-b x }
\en
for positive constants $ c $ and $b$, and $\alpha$ which is restricted
to $0 \le \alpha < 1$ by the constraints on a \Lev\ density and \re{infy}.
The corresponding distributions of $T$ are known as 
{\em generalized gamma distributions} \cite{brix99}. 
Note that the usual family of gamma distributions is recovered for
$\alpha = 0$, and that a stable distribution with index $\alpha$ is
obtained for $b = 0$ and $0 < \alpha < 1$.
One can also take $\alpha = - \kappa$ for arbitrary $\kappa > 0$,
except that in this model the \Lev\ measure has a total mass
$\psi(\infty) < \infty$ so $$\PR (T=0) = \exp(-\psi(\infty)) > 0,$$ contrary
to the present assumption that the distribution of $T$ has
a density.  Such models can be analyzed by first
conditioning on the Poisson total number of lengths,
which reduces the model to one with say
$m$ i.i.d. lengths with probability density proportional to $\rho$.
In the case \re{gamma1} for $\alpha = -\kappa$, that is to say that
the lengths are i.i.d. gamma$(\kappa,b)$ variables. This
model for random partitions has been extensively studied. 
It is well known that features of the \PDth\ model can be derived by
taking limits of this more elementary model with $m$ i.i.d.
gamma$(\kappa,b)$ lengths as $\kappa \te 0$
and $m \te \infty$ with $\kappa m \te \theta$.
See \cite{jp96bl} for a review of this circle of ideas and its applications
to species sampling models.

The \PK$(\rho_{\alpha,c,b})$ model for a random partition defined by
$\rho_{\alpha,c,b}$ in \re{gamma1} for $0 \le \alpha < 1$ was proposed by McCloskey \cite{mc65},
who first exploited the key idea of size-biased sampling in the setting
of species sampling problems.
Due to the remarks in Section \ref{secops} about scaling and exponential tilting,
for $0 < \alpha < 1$ the family of $\PK(\rho_{\alpha,c,b}, \gamma)$ distributions,
as $\gamma$ varies over all distributions on $(0,\infty)$,
depends only on $\alpha$ and not on $c$ or $b$.
So in studying this family of distributions on $\sdec$ and their associated exchangeable
partitions of $\IN$, the choice of $c$ and $b$ is
entirely a matter of convenience.
This study is taken up in the next section, with the choice of $b = 0$
and $c = \alpha / \Gamma(1-\alpha)$ which leads to the simplest form of most results.
See also \cite{brix99,james02,james02pp} regarding generalized gamma random measures and further developments.

\subsection{The stable $(\alpha)$ model}

Suppose now that $\PR_\alpha$ 
governs the Poisson model for $T$ with stable $(\alpha)$
distribution with Laplace transform
\eq
\lb{lkstable}
\ER _\alpha [ \exp ( - \lambda T) ] =  \int_0^\infty e^{-\lambda x } f_\alpha (x) \,dx =
\exp ( - \lambda ^\alpha )
\en
for some $0 < \alpha < 1 $,
where $f_\alpha(x)$ is the stable$(\alpha)$ density of $T$, 
that is \cite{po46}
\eq
\lb{aldens}
f_\alpha(t) = { - 1 \over \pi } \sum_{k = 0}^\infty  { ( - 1)^{k} \over k! } \sin ( \pi \alpha k ) { \Gamma ( \alpha k + 1 ) \over t ^{ \alpha k + 1 } }. 
\en
For $\alpha = \hf$ this reduces to the following
formula of Doetsch \cite[pp. 401-402]{doetsch37} and \Lev\ \cite{lev39}:
\eq
\lb{halfdens}
\PR_{\scrhf} ( 2 T \in dx)/dx = \hf f_{\scrhf}(\hf x ) = { 1 \over \sqrt { 2 \pi } } x ^ {-{3 \over 2}} e ^ {-{1 \over 2x}} .
\en
Special results for $\alpha = \hf$, discussed in Section \ref{sechf},
involve cancellations due to simplification of 
$f_\alpha (pt) / f_\alpha (t)$ for $0 < p < 1$,
which does not appear to be possible for general $\alpha$.
The \Lev\ density corresponding to the Laplace transform \re{lkstable}
is well known to be
\eq
\lb{lkrho}
\rho_\alpha (x ) = {  \alpha \, x^{ - \alpha - 1 } \over \Gamma(1 - \alpha) } 
~~~~~(x > 0 ).
\en
Write $\PR_ \alpha ( \cspdot |\, t )$ for 
$\PR_ \alpha ( \cspdot |\, T = t )$. 
So the $\PR_ \alpha$ distribution of $(P_i)$ on $\sdec$ is \PK$(\rho_\alpha)$,
and the $\PR_ \alpha ( \cspdot |\, t )$ distribution
of $(P_i)$ is \PK$(\rho_\alpha \giv t)$.
Note from \re{lkstable} that if $T_c$ is the total length in the model
governed by $c \rho_\alpha$ for a constant $c > 0$, then
$T_c$ has the same distribution as $c^{1/\alpha}T_1$ for $T_1 = T$ as in \re{lkstable}. Together with similar
scaling properties of the lengths $J_i$, this implies that for
all $0 < \alpha < 1$ and $t > 0$ there is the formula
\eq
\lb{scaling}
\PK(c \rho_\alpha \giv t)  = \PK(\rho_\alpha \giv c^{-1/\alpha } \,t).
\en
Formulas for the \PK$(\rho_\alpha \giv t)$ distribution are 
described in Section \ref{seccon}. 
These formulas can be understood as disintegrations 
of simpler formulas obtained in \cite{jp.epe}, and recalled in 
Section \ref{twoparam}, for a particular subfamily of the class of
 \PK $( \rho_\alpha, \gG )$ distributions.

One reason for special interest in the Kingman family associated with the
stable \Lev\ densities $\rho_\alpha$ is the following result which 
will be proved elsewhere. 

\begin{thm}{charc}
The EPPF of an exchangeable random partition $\Pi$ of $\IN$ with
an infinite number of classes with proper frequencies has an EPPF of the 
Gibbs form
\eq
\lb{gibbs}
p( n_1, \cdots, n_k ) = c_{n,k} \prod_{i = 1}^k w_{n_i} \mbox{   where $n = \sum_{i = 1}^k n_i$ }
\en
for some positive weights $w_1=1,w_2, w_3, \ldots$ 
and some $c_{n,k}$ if and only if 
$$
w_m = \prod_{j = 1}^{m-1} (j-\alpha)     ~~~~ (m = 1,2, \ldots)
$$
for some $0 \le \alpha < 1$.
If $\alpha = 0$ then the distribution of $\Pi$ corresponds to
$\int_0^\infty \PDth \gG (d \theta)$ for some probability distribution
$\gamma$ on $(0,\infty)$, whereas if $0 < \alpha < 1$ then
the distribution of $\Pi$ corresponds to 
\PK $( \rho_\alpha, \gG ):= \int_0^\infty \PK ( \rho_\alpha \giv t ) \gG (d t)$
for some $\gG$.
\end{thm}

See also Kerov \cite{kerov95} and Zabell \cite{za94} for related 
characterizations of the two-parameter family discussed
in Section \ref{twoparam}. This family is characterized by an EPPF of the
form \re{gibbs} with $c_{n,k}$ a product of a function of 
$n$ and a function of $k$.

\subsection{Conditioning on $T$}
\label{seccon}
Assume throughout this section that $0 < \alpha < 1$.
Immediately from \re{strucha} and \re{lkrho}, 
in the \PK$(\rho_\alpha \giv t)$ model, 
the distribution of $\tilP_1$ has density
\eq
\lb{lk13}
\tilde{f}_\alpha(p\giv t ) =
{\alpha ( p t )^{- \alpha }
\over \Gamma(1 - \alpha) }
\,
{f_\alpha ( (1-p) t ) \over f_\alpha ( t ) }   ~~~~~~(0 < p < 1).
\en
Let $h$ be a non-negative measurable function with $\ER_ \alpha h (T) = \int_0^\infty h(t) f_\alpha(t) dt = 1$, 
and let $h\cdot f_\alpha$ denote the distribution on $(0,\infty)$ with density $h(t) f_\alpha(t)$.
Then by integration from \re{lk13}, under the probability $\PR_ {\alpha, h}$ governing the
\PK$(\rho_\alpha , h \cdot f_\alpha )$ model, the structural distribution of $\tilP_1$ has density
\eq
\lb{lk13a}
\PR_ {\alpha,h} ( \tilP_1 \in dp)/dp = { \alpha \over \Gamma(1-\alpha)} \, p^{- \alpha } ( 1 - p )^{ \alpha - 1} \eta_{\alpha,h} ( 1-p)  ~~~~( 0 < p < 1)
\en
where
\eq
\lb{sp.31}
\eta_{\alpha,h} ( u ) := \int_0^{\infty} v^{- \alpha} h( v/u) f_{\alpha}(v) dv = \ER_ {\alpha} [ T^ {-\alpha} h ( T / u )] .
\en
For instance, it is known \cite{ppy92} that
\eq
\lb{thmom}
C_\alth := \ER_ {\alpha} ( T^{-\theta} )
= { \Gamma (\thoval + 1 ) \over \Gamma ( \theta + 1 )} ~~~~~~~~~~~~~(\theta > - \alpha).
\en
So for $\theta > - \alpha$, \re{lk13a} and \re{sp.31} imply:
\eq
\lb{lmmppy3}
\mbox{if $h(t) = C_\alth^{-1} t^{- \theta}$ then $\tilP_1$ has beta$(1-\alpha, \alpha + \theta)$ distribution.}
\en
This example is discussed further in the next section.
As another example, if $h(t) = \exp ( b^{\alpha} - b t)$ for some $b >0$, then according to
\re{pktilt} the model \PK$(\rho_\alpha , h \cdot f_\alpha )$ 
is identical to the unconditioned generalized gamma model 
\PK$(\rho_{\alpha,b})$ with 
$$
\rho_{\alpha,b}(x) := \rho_\alpha(x) e^{-bx} = 
{  \alpha \over \Gamma(1 - \alpha) } { e^{-b x }  \over x^ {\alpha + 1 } } ~~~~~(x > 0 ).
$$
So the structural density of the \PK$(\rho_{\alpha,b})$ model
is given by formula \re{lk13a} with
\eq
\lb{sp.32}
\eta_{\alpha,h}  ( u ) = \exp(b ^{\alpha}) 
\ER _\alpha [ T^ {-\alpha} \exp ( - b T / u )] .
\en
For $\alpha = \hf$ the expectation in 
(\ref{sp.32}) can be evaluated by 
using \re{halfdens} to write for $\xi >0$
\eq
\lb{sp.32a}
\ER _{\scrhf} [ T^ {- \scrhf} \exp ( - \xi T )]  =
{1 \over \sqrt{ \pi} } \int_0^{\infty} {dx \over x^2} \, e^ { - ( \xi x + 1/x)/2 }
= 2 \sqrt{ \xi \over \pi } K_1( \sqrt {\xi} )
\en
where $K_1$ is the usual modified Bessel function.
Thus for $b > 0$ the \PK$(\rho_{\scrhf,b})$ model associated with the inverse Gaussian distribution \cite{seshadri93} has structural distribution with
density $\tilde{f}_{ \scrhf,b}$ given by the formula
\eq
\lb{sp.33}
\tilde{f}_{\scrhf,b}(p) =
{ \sqrt{b} e ^{\sqrt{b} }
\over
\pi \sqrt{p} ( 1 - p ) }
K_1 \left( \sqrt { b \over  (1-p) } \right) ~~~~~~(0 < p < 1).
\en

\begin{prp}{prpstable}
For $0 < \alpha < 1, q > 0$ let
$\mu_\alpha ( q \giv t)$ denote the $q$th moment of the
structural density {\em \re{lk13}} of the \PK$(\rho_\alpha \giv t)$ 
distribution:
\eq
\lb{lk14}
\mu_\alpha ( q \giv t) 
:= \int_0 ^1 p^q  \tilde{f}_\alpha(p\giv t ) \,dp
=  \ER_  \alpha ( \tilP_1 ^q \giv t )  .
\en
Then for each $t > 0$ the EPPF of a \PK$(\rho_\alpha \giv t)$ partition of $\IN$
is
\eq
\lb{sp.19a}
p_{\alpha} \nnnt = 
\,
{\Gamma ( 1 - \alpha )  \over  \Gamma(n - k \alpha )}
\left( { \alpha \over t^\alpha }\right)^{k-1} 
\mu_\alpha ( n-1 - k\alpha + \alpha \giv t )
\,
\prod_{i= 1}^k [1 - \alpha]_{n_i - 1}
\en
where 
$$
[1 - \alpha]_{n_i-1}:= \prod_{j = 1}^{n_i-1} ( j - \alpha)  = {\Gamma (  n_i - \alpha ) \over \Gamma (1 - \alpha ) } .
$$
Alternatively,
\eq
\lb{sp.19ab}
p_{\alpha} \nnnt 
= {\alpha^k \over t^n} g_\alpha( n - k \alpha \giv t )
\prod_{i= 1}^k [1 - \alpha]_{n_i - 1}
\en
where $g_\alpha ( q \giv t ) := ( \Gamma(q) f_\alpha(t) )^{-1} \int_0^t f_\alpha(t - v ) v^{q-1} dv$.
\end{prp}
\proof
This is read from Theorem \ref{prp33}, since the integral \re{sp.18b} reduces 
to a standard Dirichlet integral.
\endpf

As checks on \re{sp.19a}, the symmetry in $\nnnn$ is still
evident, and $p_\alpha (n \giv t ) = \mu_\alpha ( n - 1 \giv t)$ as required by \re{cons4}.
However, the addition rules \re{addr} for this EPPF are not at all obvious.
Rather, they amount to the following
identity involving moments of the structural distribution:

\begin{crl}{momident}
The moments $\mu_\alpha(q | t )$ of the structural distribution
on $(0,1)$ associated with the \PK$(\rho_\alpha \giv t)$ distribution 
on $\sdec$ satisfy the following identity:
for all $1 \le k \le n$ and $ t > 0$
\eq
\lb{momalph}
\mu_\alpha ( n - 1 - k \alpha + \alpha \giv t ) =
\mu_\alpha ( n -  k \alpha + \alpha \giv t ) 
+ { \Gamma ( n - k \alpha ) \, \alpha \, t^{- \alpha} \over \Gamma (n+1 - k \alpha - \alpha ) } \,
\mu_\alpha ( n -  k \alpha \giv t )  .
\en
\end{crl}
To illustrate, 
according to the simplest addition rule \re{cons1}, 
$$
1= p_\alpha (2 \giv t ) + p_\alpha (1,1 \giv t )  ,
$$
which amounts to \re{momalph} for $n=k=1$, that is
\eq
1 = \mu_\alpha ( 1 \giv t ) + 
\, { \Gamma ( 1 - \alpha ) \over \Gamma ( 2 - 2 \alpha ) } \, 
{ \alpha  \over t ^{\alpha }}\,
\mu_\alpha ( 1 - \alpha \giv t ) .
\en
The addition rule underlying \re{momalph}
can be checked for general $\alpha$ by
an argument described in Section \ref{twoparam}.
In the case $\alpha = \hf$, 
the later formulae \re{momhf1} and \re{bmoms} show that
\re{momalph} reduces to a known recursion \re{recurh} for the
Hermite function.

Repeated application of \re{momalph} shows that for each $1 \le k \le n $
the moment on the left side of \re{sp.19a}
can be expressed as a linear combination of integer moments $\mu_\alpha (j \giv t)$ for $j = 0 , \cdots , n-1$,
with coefficients depending on $n,k,\alpha,t$
which could easily be computed recursively.
But except in the special case $\alpha = \hf$ discussed in
Section \ref{sechf}, even the integer moments seem difficult to evaluate.

\section{The two-parameter Poisson-Dirichlet family}
\label{twoparam}

For $0 < \alpha < 1, \theta > - \alpha$, let $\gG_\alth$ denote the distribution
on $(0,\infty)$
with density $C_\alth ^{-1} t^{-\theta}$
at $t$ relative to the stable$(\alpha)$ distribution of $T$ defined by \re{lkstable}, that is
\eq
\lb{falth}
\gG_\alth ( d t ) = C_\alth ^{-1} \,  t ^{-\theta } \, f_\alpha (t ) \,dt
\en
where 
$C_\alth := \ER_\alpha(T^{- \theta}) = \Gamma (\thoval + 1 )/\Gamma ( \theta + 1 )$ as in \re{thmom} and \re{lmmppy3}.
\begin{dfn}{pd2} 
{\em \cite{ppy92,py95pd2}
The {\em Poisson-Dirichlet distribution
with two parameters $(\alpha, \theta)$}, denoted \PDalth,
is the distribution on $\sdec$ defined for 
$0 \le \alpha < 1, \theta > -\alpha$ by
\eq
\mbox{ \PDalth } 
= \left\{ \begin{array}{ll}
\PDth & \mbox{ for } \alpha = 0 , \theta > 0
\\
\PK ( \rho_\alpha, \gG_\alth ) & \mbox{ for } 0 < \alpha < 1, \theta > - \alpha
\end{array}
\right.
\en
}
\end{dfn}
This family of distributions on $\sdec$ has some remarkable properties and 
applications.
As shown in \cite{ppy92}, it follows from Lemma \ref{ppy2} that
if $(P_i)$ has \PDalth distribution 
then the corresponding size-biased sequence
$(\tilP_j)$ can be represented as
\eq
\lb{sbrk}
\tilP_j = \RW_j \prod_{i=1}^{j-1} (1 - \RW_i ) 
\en
\eq
\lb{bstick}
\mbox{where the $\RW_j$ are independent with beta$(1-\alpha,\theta + j \alpha)$ distributions.}
\en
So the $\PDalth$ distribution can 
just as well be defined,
without reference to the \PKK\ construction,
as the distribution of $(P_i)$ defined by ranking $(\tilP_j)$ constructed by \re{sbrk} 
from independent $\RW_j$ as in \re{sbrk}.
The sequence $(\tilP_j)$ defined by \re{sbrk} and \re{bstick} for 
$0 \le \alpha < 1$ and $\theta > 0$ was considered by Engen \cite{en78} as a model for 
species abundances. 
See \cite{py95pd2} for further study of the \PDalth\ family.
It was shown in \cite{jp.isbp} that if $(P_i)$ is a random element of $\sdec$
with $P_i > 0 $ a.s. for all $i$ and the corresponding size-biased sequence
$(\tilP_j)$ admits the representation \re{sbrk} with
independent residual fractions $\RW_j$, then the $\RW_j$ must have beta 
distributions as described in \re{bstick},
and hence the distribution of $(P_i)$ must be
\PDalth\ for some $0 \le \alpha < 1$ and $\theta > - \alpha$.
Reformulated in terms of random partitions, 
and combined with Proposition \ref{prptk}, this yields the following:

\begin{prp}{indpceprop}
Let $\Pi$ be the exchangable random partition of $\IN$ derived by sampling from
a random element $(P_i)$ of $\sdec$ with $P_i > 0$ for all $i$.
Let $\Pi_k$ be derived from $\Pi$ by deletion of the first $k$ classes of $\Pi$,
with classes in order of appearance, as defined above Proposition 
{\em \ref{prptk}}.
Then the following are equivalent

\noindent{\em (i)}
for each $k$, $\Pi_k$ is independent of the frequencies
$(\tilP_1, \cdots, \tilP_k)$ of the first $k$ classes of $\Pi$;

\noindent{\em (ii)}
$\Pi$ is a \PDalth-partition  for some 
$0 \le \alpha < 1$ and $\theta > - \alpha$, in which case 
$\Pi_k$ is a $\PD(\alpha, \theta + k \alpha )$-partition.
\end{prp}

As shown in \cite{jp.epe},
the independence property \re{sbrk} of the residual fractions
$\RW_j$ of a \PDalth-partition allows the corresponding EPPF 
$p_\alth(n_1 , \ldots , n_k )$
to be evaluated using \re{genepf}.
The result is as follows.
For all $0 \le \alpha <1$
and $\theta > - \alpha$, 
\eq
\lb{class.esf}
p_\alth(n_1 , \ldots , n_k )
= { [\theta + \alpha ]_{k-1;\alpha} \over 
[ \theta + 1 ]_{n-1 } }
 \prod_{i=1}^k [1- \alpha]_{n_i - 1}
\en
where $n = \sum_{i = 1}^k n_i$ and
for real $x$ and $a$ and non-negative integer $m$
$$
[ x ]_ {m;a} 
= \left\{ \begin{array}{l}
1~\mbox{ for }~m=0 \\
x (x + a ) \cdots (x +(m-1) a )~\mbox{ for }~m=1,2, \ldots
\end{array}
\right.
$$
and $[x]_m = [x]_{m;1}$.
The previous formula \re{esf} is the special case of
\re{class.esf} for $\alpha = 0$. 
Both this case of \re{class.esf}, and the case 
when $0 < \alpha < 1$ and $\theta = 0$, follow easily from \re{eppf2}.
Formula
\re{class.esf} shows that a \PD$(\alpha, \theta )$ partition $\Pi$
of $\IN$ to be constructed sequentially as follows \cite{jp.epe,jp96bl}. 
Starting from $\Pi_1 = \{ \{1 \} \}$, given that $\Pi_n$ has been
constructed as a partition of $\IN_n$ with say $k$ blocks of sizes
$(n_1, \cdots, n_k)$, define $\Pi_{n+1}$ by assigning the new element
$n+1$ to the $j$th class whose current size is $n_j$ with
probability
\eq
\lb{pred1}
\PR ( j \uparrow \giv n_1, \cdots , n_k ) = { n_j - \alpha  \over 
 n + \theta }
\en
for $1 \le j \le k$,
and assigning $n+1$ to a new class numbered $k+1$ with the remaining probability
\eq
\lb{pred2}
\PR ( k+1 \uparrow \giv n_1, \cdots , n_k ) = { k \alpha  \over  n + \theta  }
\en
For $\alpha = 0$ and $\theta > 0$ this is generalization of Polya's
urn scheme developed by Blackwell-McQueen \cite{bl73} and Hoppe \cite{ho87}. 
See \cite{jp.epe,jp96bl,hansenp98} for consideration of more general
prediction rules for exchangeable random partitions.

The following calculation shows how to derive either of the
two EPPF's \re{class.esf} and \re{sp.19a} from the other.
The argument also shows that the function
$\paltnnn$ defined by \re{sp.19a} 
satisfies the addition rules of an EPPF 
as a consequence of the corresponding addition rules for $\palthnnn$,
which are much more obvious.

The kernel
$\gG_\alth (dt) $ introduced in \re{falth}, is now viewed for
a fixed $\alpha$ as a family of probability distributions
on $(0,\infty)$ indexed by $\theta \in (-\alpha, \infty)$,
that is a Markov kernel $\gG_\alpha$ from $(-\alpha, \infty)$ to 
$(0,\infty)$.
For a non-negative measurable function 
$h = h(t)$ with domain $(0,\infty)$, define
a function $\gG_\alpha h = (\gG_\alpha h )(\theta )$ with domain $(-\alpha, \infty)$ by the usual action of this Markov kernel as an integral operator:
\eq
(\gamma_\alpha h ) (\theta) =
\int_0^\infty \gG_\alth (dt) h(t)
\en
Then say $(\gamma_\alpha h )(\theta) $ is the {\em $\gamma_\alpha$-transform } of $h(t)$.
Let $\ER_ \alth$ denote expectation with respect to the
probability distribution 
$$\PR_ \alth (\cdot) := \int_0 ^ \infty \PR_ \alpha (\cdot \giv t ) \gG_\alth(dt).$$
By definition, for each non-negative random variable $X$ governed by
the family of conditional laws $(\PR_ \alpha (\cdot \giv t ), t >0)$,
\eq
\mbox{ the $\gamma_\alpha$-transform of } 
\ER_ \alpha( X \giv t ) 
\mbox { is } 
\ER_ {\alth}(X ).
\en
In particular, for each $\nnnn$,
\eq
\mbox{ the $\gamma_\alpha$-transform of } 
\paltnnn
\mbox { is } 
\palthnnn .
\en
An obvious change of variable allows uniqueness and inversion
results for the $\gamma_\alpha$-transform to be deduced from standard results
for Mellin or bilateral exponential transforms.
So the problem is just to show that
the $\gamma_\alpha$-transform of the right side of \re{sp.19a} is the right side of \re{class.esf}.
%
%
%
%
To see this, observe first that for each $q > 0$, because
$\mu_\alpha (q | t )  := \ER_ \alpha ( \tilP_1 ^q \giv t)$,
\eq
\lb{betatf}
\mbox{ the $\gamma_\alpha$-transform of } 
\mu_\alpha (q | t )  
\mbox{ is  }
\ER_ \alth ( \tilP_1 ^q ) = {\Gamma ( 1 - \alpha  + q ) \Gamma ( 1 + \theta )
\over
\Gamma ( 1 + \theta  + q ) \Gamma ( 1 - \alpha) }
\en
where $\ER_ \alth ( \tilP_1 ^q )$ is evaluated using \re{lmmppy3}.
To deal with the factor of $t^{- (k-1) \alpha }$ in \re{sp.19a},
note from \re{thmom} that for each $\beta > 0$, and any $h(t)$,
\eq
\lb{hal88}
\mbox{the $\gamma_\alpha$-transform of $t^{-\beta } h(t)$ is  }
{\Gamma( \thoval + \beoval + 1 ) \Gamma ( \theta + 1 )
\over
\Gamma( \thoval + 1 ) \Gamma ( \theta + \beta + 1 ) } \, 
\, 
(\gamma_\alpha h )(\theta + \beta)  .
\en
By \re{betatf} for $q = n -1 - k \alpha + \alpha $ and 
\re{hal88} for $\beta =  \alpha k - \alpha $ and 
$h(t) = \mu_\alpha(q \giv t )$
the right side of \re{sp.19a}
has for its $\gamma_\alpha$-transform the following function of $\theta$:
$$
{ \alpha ^{k-1} \Gamma(1-\alpha)
\over
\Gamma ( n - k \alpha ) }
\,
{\Gamma( \thoval + k ) \Gamma ( \theta + 1 )
\, 
\over
\Gamma( \thoval + 1 ) \Gamma ( \theta + k\alpha - \alpha +  1 ) } 
\, 
\,
{ \Gamma  (n - k \alpha  ) \Gamma ( 1 + \theta + k \alpha - \alpha )
\over
\Gamma ( n + \theta ) \Gamma ( 1 - \alpha ) }
\, 
\prod_{i= 1}^k [1 - \alpha]_{n_i - 1}
$$
which reduces by cancellation to the right side of \re{class.esf}.

\subsection{The $\alpha$-diversity}

Let $\Pi$ be an exchangeable random partition of $\IN$ with ranked frequencies
$(P_i)$. Let $K_n$ denote the number of classes of $\Pi_n$,
the partition of $\IN_n$ induced by $\Pi$.
Say that $\Pi$ has {\em $\alpha$-diversity $S$} and write
$\aldiv ( \Pi) = S$
iff there exists a random variable $S$ with 
$0 < S < \infty $ a.s. and
\eq
K_n  \sim  S n ^{\alpha}  \mbox{ as } n \te \infty
\en
where for two sequences of random variables $A_n$ and $B_n$, the notation
$A_n \sim B_n$ will now be used to indicate that $A_n/B_n \te 1$ almost
surely as $n \te \infty$.
According to a result of Karlin \cite{kar67urn},
applied conditionally given $(P_i)$,
if these ranked frequencies are such that 
\eq
\lb{poi.4bb}
P_{i}
\sim  \left( {S \over  \Gamma ( 1 - \alpha ) i } \right) ^ { {1 \over \alpha } }  
\en
for some $0 < S < \infty$
then $\Pi$ has $\alpha$-diversity $S$. 

\begin{prp}{aldiv}
Suppose $\Pi$ is a $\PK(\rho_\alpha, \gG )$ partition of $\IN$ for some
$0 < \alpha < 1$ and some probability distribution $\gG$ on $(0,\infty)$.
Then

\noindent
\rm{(i)}
$\aldiv ( \Pi) = S$ for a random variable $S$ with $S = T^{-\alpha}$
where $T = S^{-1/\alpha}$ has distribution $\gG$. In particular, $S = t^{-\alpha}$ is constant if $\Pi$ is a $\PK(\rho_\alpha \giv t )$ partition.

\noindent
{\rm (ii)} A regular conditional distribution for $\Pi$ given $S = s$ is defined 
by the EPPF $p_\alpha ( n_1, \cdots, n_k |s^{-1/\alpha} )$ 
obtained by setting $t = s^{-1/\alpha}$ in \re{sp.19a}.

\noindent
{\rm (iii)}
In particular, both {\rm (i)} and {\rm (ii)}
hold if $\Pi$ is a $\PD(\alth)$ partition for some $\theta > - \alpha$.
Then the $\alpha$-diversity $S$ of $\Pi$ is $S = T^{-\alpha}$ for $T$ with
the distribution $\gG_\alth $ defined by \re{falth}.
\end{prp}
\proof
Suppose that $(P_i)$ has $\PK(\rho_\alpha, \gG )$ distribution.
The fact that \re{poi.4bb} holds for $S = T^{-\alpha}$ in the
unconditioned case where $T$ has stable$(\alpha)$ distribution is due to
Kingman \cite{ki75}. Kingman's argument, using the law of large numbers
for small jumps of the Poisson process, applies just as well for $T$
conditioned to be a constant $t$. So \re{poi.4bb} follows in general
by mixing over $t$. 
\endpf

See \cite{py95pd2} and papers cited there for further information about the 
Mittag-Leffler distribution of $S = T^{-\alpha}$ derived from a $\PD(\alpha, 0)$
partition. The corresponding distribution of $S$ for $\PD(\alth)$ has density 
at $s$ proportional to $s^{\theta/\alpha}$ relative to this Mittag-Leffler distribution. 

As shown in \cite[Proposition 10]{py95pd2},
if $\Pi$ is a partition of $\IN$ whose ranked frequencies $(P_i)$
have the $\PD(\alpha,0)$ distribution, then $S = \aldiv ( \Pi) $ can be recovered from 
$\Pi$ or $(P_i)$ via either \re{hal88} or \re{poi.4bb}. Then $T = S^{-1/\alpha}$ 
has stable$(\alpha)$ distribution as in \re{lkstable}, and 
$(TP_i)$ is then sequence of points of a Poisson process with \Lev\
density $\rho_\alpha$. 
See also \cite{jp99sam,csp} for more about the distribution of $K_n$
derived from a $\PD(\alpha,\theta)$ partition.

\section{Application to lengths of excursions.} 
\label{excns}
This section reviews some results of \cite{ppy92,py92,jp.bmpart,py95pd2}.
Let $\PO$ govern a strong Markov process $B$ starting at a recurrent point 0
of its statespace, 
such that the inverse $(\tau_\ell, \ell \ge 0 )$
of the local time process $(L_t, t\ge 0 )$ of $B$ at zero 
is a stable subordinator of index $\alpha$ for some $0 < \alpha < 1$.
That is to say,
$\EO \exp ( - \lambda \tau_1 ) = \exp ( - c \lambda^\alpha )$ for
some constant $c > 0 $.
So the $\PO$ distribution of $\tau_1$ is the $\PR_ \alpha$ distribution
of $c^{1/\alpha}  T $ for $T$ as in \re{lkstable}.
For example,
$B$ could be a one-dimensional Brownian motion $(\alpha = \hf)$ or Bessel process of dimension $2 - 2 \alpha$.
In the Brownian case, take $c = \sqrt{2}$ to obtain
the usual normalization of local time as occupation density relative to Lebesgue measure,
which makes $L_1 \ed |B_1|$.
Let $M = \{t : 0 \le t \le 1 , B_t = 0 \}$ denote
the random closed subset of $[0,1]$ defined by the zero set of $B$.
Component intervals of the complement of $M$ relative to $[0,1]$ are
called {\em excursion intervals}.
For $0 \le t \le 1 $ let $G_t = \sup\{M \cap [0,t]\}$, the last
zero of $B$ before time $t$.
Note that with probability one, $G_1 < 1$, so one of the
excursion intervals is the {\em meander interval} $(G_1,1]$,
whose length $1-G_1$ is one of the lengths appearing in the list
$(P_i)$ say of ranked lengths of excursion intervals.
According to the main result of \cite{py92},
\eq
\lb{ex1}
\mbox{the sequence $(P_i)$ of ranked lengths has $\PD (\alpha,0)$ distribution}
\en
Let $U_1, U_2, \cdots  $ be a sequence of i.i.d. uniform $[0,1]$
random variables, independent of $B$, called the sequence of
{\em sample points}.
Let $\PIB{}  = (\PIB{n})$
be the random partition of $\IN$ generated by the random equivalence
relation $i \sim j $ iff $G_{U_i} = G_{U_j}$.
That is to say $i \sim j$ iff $U_i$ and $U_j$ fall in the same excursion 
interval. So for example $\PIB{5} = \{ \{1,2,5\} , \{3\}, \{4\} \}$
iff $U_1, U_2 $ and $U_5$ fall in one excursion interval,
$U_3$ in another, and $U_4$ in a third.
By translation of results of \cite{py92,py95pd2} into present notation
\eq
\lb{ex2}
\mbox{ $\PIB{}$ is a $\PD (\alpha,0)$ partition and $\aldiv ( \PIB {}) = c L_1$} 
\en
where $L_1$ is the local time of $B$ at zero up to time 1.
By construction, the sequence $(\tilP_j)$ of class frequencies of $\PIB{}$
is the sequence of lengths of excursion intervals
in the order they are discovered by the sample points,
and $(P_i)$ is recovered from $(\tilP_j)$ by ranking.
To illustrate formula \re{class.esf}, $U_1$ and $U_2$ fall in different excursion
intervals with probability $p_{\alpha,0}(1,1) = \alpha$, and in the same one
with probability $p_{\alpha,0}(2) = 1-\alpha$.
%
Similarly,
given that the local time is $L_1  = \ells$,
two sample points fall in the same excursion interval
with probability $p_\alpha(2 \giv ( c \ells) ^{-1/\alpha})$, and
in different excursion intervals with probability
$p_\alpha(1,1 \giv (c \ells)^{-1/\alpha})$, for $p_\alpha( \cdots \giv t )$
defined by \re{sp.19a}.
See Section \ref{sechf} for evaluation of these functions in the case $\alpha = \hf$
corresponding to a Brownian motion $B$.

Let $\tilR_n = 1 - \tilP_1 - \cdots - \tilP_{n}$, which is the total length of
excursions which remain undiscovered after the sampling process has found $n$
distinct excursion intervals.
The result of Proposition \ref{indpceprop} in this setting, due to \cite{ppy92},
is that for each $n = 0, 1, 2, \cdots$ a $\PD(\alpha, n \alpha )$ distributed
sequence is obtained by ranking the sequence 
\eq
\lb{ex3}
{1 \over \tilR_{n} } ( \tilP_{n+1} , \tilP_{n+2},  \cdots )
\en
of relative excursion lengths which remain
after discovery of the first $n$ intervals.
For $n=1$ the same $\PD(\alpha, \alpha )$ distribution is
obtained more simply by deleting the meander of length $1-G_1$,
renormalizing and reranking. 
This is due to the result of \cite{py92} that the length $1-G_1$ of
the meander interval is a size-biased choice from $(P_i)$.
As the excursion lengths in this case are just
the excursion lengths of a standard bridge, equivalent to conditioning on 
$B_1 = 0$, the ranked excursion lengths of such a bridge have 
$\PD(\alpha,\alpha)$ distribution. As first shown in \cite{py92},
this implies that both the unconditioned process $B$
and the bridge $B$ given $B_1 = 0$ share a 
common conditional distribution for the ranked excursion lengths
$(P_i)$ given the local time $L_1$.
In present notation, this conditional distribution of $(P_i)$
given $L_1 = \ell$, with or without conditioning on $B_1 = 0$,
is $\PK(\rho_\alpha | (c \ells)^{-1/\alpha})$.


One final identity is worth noting. As a consequence of the above discussion,
for the process $B$, the conditional distribution of 
the meander length $1-G_1$ given $L_1 = \ell$ is given by
\eq
\PO ( 1 - G_1 \in dp | L_1 = \ell ) =  \PO ( \tilP_1 \in  \ dp | L_1 = \ell ) = 
\tilde{f}_\alpha ( p | (c \ell )^{-1/\alpha} ) dp
\en
where $\tilde{f}_\alpha ( p | t)$
as in \re{lk13} is the structural density 
of the Poisson model for stable $(\alpha )$ distributed $T$ 
conditioned on $T = t$.
So the moment function $\mu_\alpha(q \giv t )$ appearing in the EPPF \re{sp.19a}
of this model can be interpreted in the present setting as
\eq
\mu_\alpha ( q \giv t ) = \EO [ (1 - G_1)^q \giv L_1 = c^{-1} t ^{-\alpha}].
\en

\section{The Brownian excursion partition}
\label{sechf}
In this section let 
$\Pi$ be the {\em Brownian excursion partition},
that is the random partition of $\IN$ 
generated by uniform random sampling of points from the interval
$[0,1]$ partitioned by the excursion intervals of  a standard Brownian 
motion $B$. According to the result of \cite{py92} recalled in \re{ex1},

\eq
\lb{pkhf1}
\Pi \mbox{ is a $\PK(\rho_{\scrhf}) = \PD(\hf,0)$ partition.}
\en
With conditioning on $B_1 = 0$, the process $B$ becomes a
standard Brownian bridge. So $\Pi$ given
$B_1 = 0$ is a $\PD(\hf,\hf)$ partition, as discussed in the previous subsection.
Features of the distribution of $\Pi$ and the conditional
distribution of $\Pi$ given $B_1 = 0$
were described in \cite{jp.bmpart}. This section presents refinements of these
results  obtained by conditioning on $L_1$, the local time of $B$
at $0$ up to time $1$, with the usual normalization of Brownian
local time as occupation density relative to Lebesgue measure.
Unconditionally, $L_1$ has the same distribution as $|B_1|$, that is
$$
\PR ( L_1 \in d \ll ) = \PR ( |B_1| \in d \ll ) = 2 \varphi(\ll) d \ll ~~~~~~~~~~~~~~ (\ll >0)
$$
where $\varphi (z ) := ( 1/\sqrt{ 2 \pi }) \exp(-\hf z^2)$ is the standard
Gaussian density of $B_1$.
Whereas the conditional distribution of $L_1$ given $B_1 = 0$ is
the Rayleigh distribution
$$
\PR ( L_1 \in d \ll \giv B_1 = 0 )  = \sqrt{ 2 \pi }  \ll \varphi(\ll) d \ll 
~~~~~~~~~~~~ (\ll >0) .
$$
Note from \re{ex2} that the $\hf$-diversity of $\Pi$ is the random variable $\sqrt{2}L_1$.
So the number $K_n$ of blocks of $\Pi$ grows almost surely like 
$\sqrt{2 n} L_1$ as $n \te \infty$.
For $\ll \ge 0$ let $\Pi(\ll)$ denote a random partition with
\eq
\lb{eqdis}
\Pi(\ll) \ed (\Pi \giv L_1 = \ll ) \ed (\Pi \giv L_1 = \ll , B_1 = 0)
\en
where $\ed$ denotes equality in distribution. 
So according to the previous discussion, 
\eq
\lb{pkhf}
\Pi(\ll) \mbox{ is a $\PK(\rho_{\scrhf} \giv  \hf \ll^{-2} )$ partition }
\en
whose $\hf$-diversity is $\sqrt{2}\ll$.
Let $\PD(\hf \bgiv \ll)$ denote the probability distribution on $\sdec$ associated
with $\Pi(\ll)$, that is the common distribution of ranked lengths of excursions of a
Brownian motion or Brownian bridge over $[0,1]$ given $L_1 = \ll$. 
Then by Definition \ref{pd2} and \re{halfdens}, for $\theta > - \hf$ there is the
identity of probability laws on $\sdec$
\eq
\lb{pdint}
\PD( \hf , \theta ) =  { 2 \over \ER [ |B_1|^{2 \theta}  ] } \, \int_0^\infty 
\rex {\ll} \ll^ {2 \theta }  \varphi ( \ll ) d \ll 
\en
where,
according to the gamma$(\hf)$ distribution of $\hf B_1^2$ and the duplication formula for the gamma function,
\eq
\lb{bmoms}
\ER ( |B_1|^{2 \theta} )  = 2^{\theta } { \Gamma( \theta + \hf) \over \Gamma(\hf) } = 2^{- \theta} { \Gamma( 2 \theta + 1 ) \over \Gamma(\theta + 1) }
~~~~~~~(\theta > - \hf).
\en
It was shown in \cite{jpda97sac} (see also \cite{bertoin00f,csp})
that it is possible to construct the Brownian excursion 
partitions as a partition valued {\em fragmentation process}
$(\Pi(\ll), \ll \ge 0)$,
meaning that $\Pi(\ll)$ is constructed for each $\ll$ on the same probability 
space, in such a way that $\Pi(\ll)$ is a coarser partition than $\Pi(\mu)$ whenever $\ll < \mu$.  
The question of whether a similar construction is possible
for index $\alpha$ instead of index $\hf$ remains open. A natural
guess is that such a construction might be made with one of the
self-similar fragmentation processes of Bertoin \cite{bert02ss}, but
Miermont and Schweinsberg \cite{miersch01} have recently shown that a construction of this form is possible only for $\alpha = \hf$.

\subsection{Length biased sampling}
Let $\tilP_j(\ll)$ denote the frequency of the $j$th class of $\Pi(\ll)$.
So $(\tilP_j(\ll), j = 1,2 \ldots)$ is distributed like the lengths of excursions of $B$ over $[0,1]$ given $L_1 = \ll$,
as discovered by a process of length-biased sampling.
In view of \Lev's formula \re{halfdens} for the stable$(\hf)$ density, the formula \re{lk13} reduces for
$\alpha = \hf$ to the following more explicit formula for the structural density of 
$\Pi(\ll)$:
\eq
\lb{Ydens}
\PR ( \tilP_1(\ll) \in dp) = {\ll \over \sqrt{ 2 \pi}  }
p^{ - \scrhf } (1-p)^{- \scrthf } \exp \left( - {\ll^2 \over 2 } {  p \over (1-p) } \right)  dp
~~~~~(0 < p < 1)
\en
or equivalently
\eq
\lb{Ydf}
\PR ( \tilP_1 \leq y) = 2 \Phi \left( \ll   \sqrt{ y \over 1-y } \right) - 1  
~~~~~(0 \le  y < 1) 
\en
where $\Phi(z):= \PR (B_1 \le z)$ is the standard Gaussian distribution function.
Put another way, there is the equality in distribution
\eq
\lb{Yed}
\tilP_1(\ll) \ed \frac{B_1^2}{\ll ^2 + B_1^2} .
\en
Furthermore, by a similar analysis using Lemma \ref{ppy},  there is the 
following result which shows how to construct the whole sequence
$( \tilP_j (\ll), j \ge 1 )$ for any $\ll >0$ from a single sequence of
independent standard Gaussian variables. Then $\Pi(\ll)$ can be constructed
by sampling from $( \tilP_j (\ll), j \ge 1 )$ as discussed in Section \ref{secstruc}.

\begin{prp}{prpac} {\em \cite[Corollary 5]{jpda97sac}}
Fix $\ll >0$. A sequence $( \tilP_j (\ll), j \ge 1 )$ is distributed like a length-biased
random permutation of the lengths of excursions of a Brownian motion or standard Brownian
bridge over $[0,1]$ conditioned on $L_1 = \ll$ if and only if
\eq
\lb{yst}
\tilP_j (\ll)  = {\ll ^2 \over \ll^2 + S_{j-1} } - {\ll ^2 \over \ll^2 + S_{j} }
\en
where $S_j:= \sum_{i = 1}^j X_i$ for $X_i$ which are independent and identically
distributed like $B_1^2$ for a standard Gaussian variable $B_1$.
\end{prp}
Let $\mu(q \bgiv \ll)$ denote the $q$th moment of the distribution of $\tilP_1(\ll)$.
So in the notation of \re{lk14} and \re{momalph}
\eq
\lb{momhf}
\mu(q \bgiv \ll) 
:= \ER [ (\tilP_1(\ll))^q ] = \mu_\scrhf ( q \giv  \hf \ll ^{-2} ) .
\en

\begin{lmm}{lmmq}
For each $\ll >0$ 
\eq
\lb{momhf1}
\mu(q \bgiv \ll) = \ER \left[ \left(  B_1^2 \over \ll^2 + B_1^2 \right)^q \right] =  
\ER ( |B_1|^{2q} ) \, h_{-2q} (\ll)
~~~~( q > - \hf )
\en
where $\ER ( |B_1|^{2q} ) $ is given by {\em \re{bmoms}}
and $h_{-2 q }$ is the Hermite function of index $-2 q$, that is
$h_0(\ll) = 1$ and for $q \notin \{ 0,1,2 \ldots \}$
\eq
\lb{hermq}
h_{ - 2 q } ( \ll ) := { 1 \over 2 \Gamma (2 q ) } \sum_{ j = 0}^\infty \Gamma ( q + j/2 ) 2^{ q + j/2} { (-\ll )^{j} \over j! } .
\en
Also,
\eq
\lb{gamlt2}
\mu(q \bgiv \ll) = 
\ER [ \exp ( - \ll \sqrt{  2 \Gamma_q  }  ]  ~~~~~~(q > 0 )
\en
where $\Gamma_q$ denotes a Gamma random variable with parameter $q$:
$$
\PR ( \Gamma_q \in dt ) = \Gamma(q)^{-1} t^{q-1} e^{- t } dt ~~~~~~~(t >0) .
$$
\end{lmm}
\proof
The first equality in \re{momhf1} is read from \re{Yed}. The second equality
in \re{momhf1} is the integral representation of the Hermite
function provided by Lebedev \cite[Problem 10.8.1]{lebedev65}, and \re{hermq}
is read from \cite[(10.4.3)]{lebedev65}.
According to another well known integral representation of the Hermite function
\cite[(10.5.2)]{lebedev65}, \cite[8.3 (3)]{erdelyi53II}, for $q >0$
\eq
\lb{hmp1}
h_{-2q}(x) = 
{1 \over \Gamma(2 q ) } \int_0^\infty t^{2 q -1} e^{- \hf t^2  - x t } dt 
= {2 ^{q-1} \over \Gamma(2 q ) } \int_0^\infty v^{q -1} e^{- v  - x \sqrt{ 2 v}  } d v.
\en
Formula \re{gamlt2} follows easily from this and \re{momhf1}.
\endpf

The identity
\eq
\lb{newid}
\ER \left[ \left(  B_1^2 \over \ll^2 + B_1^2 \right)^q \right] =  \ER [ \exp ( - \ll \sqrt{  2 \Gamma_q  }  ]  ~~~~(q >0),
\en
which is implied by the previous proposition, can also be checked by the following argument suggested by Marc Yor.
Let $X$ be a positive random variable independent of $\Gamma_q$, and let $\eps$ with $\eps \ed \Gamma_1$ be
a standard exponential variable independent of both $X$ and $\Gamma_q$. Then by elementary conditioning arguments,
for $\theta \ge 0$
\eq
\lb{cond1}
\ER \left[ \left( { X \over \theta + X } \right)^q \right]
= \ER \left[ e^{- \theta \Gamma_q/X } \right] = \PR ( \eps X /\Gamma_q > \theta ) .
\en
Take $X = B_1^2$ and $\theta = \ll^2$, and use the identity $\eps B_1^2 \ed \eps^2/2$, which is a well known
probabilistic expression of the gamma duplication formula, to deduce \re{newid} from \re{cond1}. 

The following display identifies $h_\nu(z)$ in the notation of various authors:
\begin{eqnarray*}
h_\nu(z) & = & 
2^{-\nu/2} H_\nu(z/\sqrt{2}) =  2^{\nu/2} \Psi( - \nu/2, 1/2, z^2/2 )
~~~~~~~~\mbox{(Lebedev\cite{lebedev65})} \\
& = & 2^{\nu /2 } U( - \nu/2, 1/2, z^2/2 )
~~~~~~~~~~~~~~~~~~~~\mbox{(Abramowitz and Stegun) \cite{as65}} \\
& = & e^{\quart z^2} U(-\nu - \hf, z) ~~~~~~~~~~~~~~~~~~~~~~~~~~~\mbox{(Miller\cite{miller55})} \\
& = & e^{\quart z^2} D_\nu(z)  ~~~~~~~~~~~~~~~~~~~~~~~~~~~~~~~~~~~~\mbox{(Erdelyi \cite{erdelyi53II}, Toscano \cite{toscano71p})} \\
\end{eqnarray*}
The functions $U(a,z)$ and $D_\nu(z)$ are known as {\em parabolic cylinder
functions}, {\em Weber functions} or {\em Whittaker functions}.
The function $U(a,b,z)$, which is 
available in {\em Mathematica} as {\tt HypergeometricU[a,b,z]},
is a {\em confluent hypergeometric function of the second kind}. 
Note that $h_n(z)$ defined for $n = 0,1,2, \ldots$ by continuous extension of \re{hermq}
is the sequence of Hermite polynomials orthogonal with respect to the standard Gaussian
density $\varphi(x)$. Also, the function $h_{-1}(x)$ for real $x$ is identified 
as {\em Mill's ratio} \cite[33.7]{JK70}:
\eq
\lb{mills}
h_{-1}(x) 
= {\PR ( B_1 > x ) \over \varphi(x) } = e^{\scrhf x^2 } \int_x^\infty e^{-\scrhf z^2 } dz.
\en
For all complex $\nu$ and $z$,
the Hermite function satisfies the recursion
\eq
\lb{recurh}
h_{\nu + 1} ( z ) = z h_\nu(z) - \nu h_{\nu - 1 } (z),
\en
which combined with \re{mills} and $h_0(x) = 1$ yields
\eq
\lb{psi2}
h_{-2} (x ) = 1 - x h_{-1}( x ) 
\en
\eq
2! h_{-3} (\xell) = - \xell + (1 + \xell^2 ) h_{-1} (\xell )
\en
\eq
3! h_{-4}(\xell) = 2+\xell ^2  - (3\xell +\xell ^3) h_{-1} (\xell )
\en
and so on.
See \cite{py02h} for further interpretations of the Hermite function
in terms of Brownian motion and related stochastic processes.

\subsection{Partition probabilities}
Recall the notation
$$
[ \hf ]_{n} :=  \prod_{j = 1}^{n} (j - \hf)  = { \Gamma ( \hf + n ) \over \Gamma ( \hf ) } = 
{ ( 2 n)!  \over 2^{2n} n ! }  .
$$

\begin{crl}{crlbm}
The distribution of $\Pi(\ll)$, a Brownian excursion partition
conditioned on $L_1 = \ll$, is determined by the following EPPF:
for $n_1, \ldots, n_k$ with $\sum_{i = 1}^k n_i = n$
\eq
\lb{eppf}
p_{\scrhf} (n_1, \ldots, n_k \bgiv \ll) = 2^{n-k} \ll ^{k-1} h_{k+1 - 2n } ( \ll ) \, \prod_{i=1}^k [ \hf ]_{n_i - 1}  .
\en
\end{crl}
\proof
This is read from \re{sp.19a}, \re{momhf1} and \re{bmoms}. 
\endpf
Formula \re{eppf} combined with \re{qn22} gives an expression in 
terms of the Hermite function for the positive integer moments of the sum 
$S_m(\la)$ of $m$th powers of lengths of excursions of Brownian motion on 
$[0,1]$ given $L_1 = \la$. This formula for $m= 2$ was derived in another
way by Janson \cite[Theorem 7.4]{janson01h}.
There the distribution of $S_2(\la)$ appears as the asymptotic distribution, 
in a suitable limit regime, of the cost of linear probing hashing. 

According to \re{pkhf} and Definition \ref{pd2},
for each $\theta > - \hf$, 
the EPPF  \re{eppf} describes the conditional distribution
of a \PD$(\hf, \theta )$ partition $(\Pi_n)$ given 
$\lim_n K_n /\sqrt{ 2 n } = \ll$, where $K_n$ is the 
number of blocks of $\Pi_n$.
Easily from \re{eppf},
for each fixed $\la > 0$, a sequential description of 
$(\Pi_n(\la), n = 1,2, \ldots)$ is obtained by replacing the prediction rules 
\re{pred1} and \re{pred2} by
\eq
\lb{pred3}
~~~~~\PR ( j \uparrow \giv n_1, \cdots , n_k) = 
( 2 n_j - 1 ) \, { h_{k -1 - 2n }(\ll)  \over h_{k + 1 -2n } (\ll) }
~~~~(1 \le j \le k)
\en
\eq
\lb{pred4}
\PR ( k+1 \uparrow  \giv n_1, \cdots , n_k) 
=
{ \ll h_{k - 2n}(\ll)  \over h_{k+ 1 -2n } (\ll) } .
 ~~~~~~~~~~~~~~~~~~~~~~~~~~~~~
\en
The addition rule for the EPPF \re{eppf} is equivalent to the
fact that these transition probabilities sum to $1$. As a check, this 
is implied the recurrence formula \re{recurh} for the Hermite function.

\begin{crl}{crlknl}
Let $K_n(\ll)$ be the number of blocks of $\Pi_n(\ll)$,
where $(\Pi_n(\ll), n = 1,2, \ldots)$ is the Brownian excursion partition
conditioned on $L_1 = \ll$.
Then 
$(K_n(\la), n = 1,2, \ldots)$ is a Markov chain with the following
inhomogeneous transition probabilities: for $1 \le k \le n$
\eq
\lb{pred33}
\PR ( K_{n+1}(\ll) = k \giv K_n(\ll) = k ) = 
( 2 n - k ) \, { h_{k -1 - 2n }(\ll)  \over h_{k + 1 -2n } (\ll) }
\en
\eq
\lb{pred44}
\PR ( K_{n+1}(\ll) = k +1 \giv K_n(\ll) = k ) = 
{ \ll h_{k - 2n}(\ll)  \over h_{k+ 1 -2n } (\ll) } .
~~~~~~~~~~~~~~~
\en
Moreover, the distribution of $K_n(\la)$ is given by the formula
\eq
\lb{fxy}
\PR (K_n (\ll) = k ) =  
\,
{(2n-k-1)! \, \ll^{k-1} h_{k+1-2n}  (\ll) \over 
(n - k) ! (k-1 )! 2^{n-k}} 
~~~~(1 \le k \le n).
\en
\end{crl}
\proof The Markov property of $(K_n(\la), n = 1,2, \ldots)$ and the
transition probabilities \re{pred33}--\re{pred44} follow easily from
\re{pred3}--\re{pred4}. 
Then \re{fxy} follows by induction on $n$, using the forwards
equations implied by the transition probabilities.
\endpf

Let $K_n$ denote the number of blocks of $\Pi_n$, where $(\Pi_n)$ is
the unconditioned Brownian excursion partition.
Then, from the discussion around \re{eqdis},
\eq
\lb{condmc}
(K_n(\la), n \ge 1 ) \ed (K_n, n \ge 1 \giv \lim_n K_n/\sqrt{2 n } = \la ).
\en
According to \re{pkhf1}, \re{pred1} and \re{pred2}, 
the sequence
$(K_n, n \ge 1)$ is an inhomogeneous  Markov chain with transition probabilities
\eq
\lb{pred333}
\PR ( K_{n+1} = k \giv K_n = k ) =  {2 n - k \over 2 n }
\en
\eq
\lb{pred444}
\PR ( K_{n+1} = k +1 \giv K_n = k ) =  { k \over 2 n } ~~~~~~~~
\en
which imply that the unconditional distribution of $K_n$ 
is given by the formula \cite[Corollary 3]{jp.bmpart} 
\eq
\lb{brownkn}
\PR (K_n = k ) = { 2n - k - 1 \choose n - 1 } 2 ^{ k+1 - 2n } 
~~~~(1 \le k \le n).
\en
Due to \re{condmc}, for each $\la > 0$ the inhomogeneous Markov chain 
$(K_n(\la), n \ge 1)$ has the same co-transition probabilities as
$(K_n, n \ge 1)$. From \re{pred333}, \re{pred444} and \re{brownkn},
the co-transition probabilities of $(K_n, n \ge 1)$ are
\eq
\lb{pred3333}
\PR ( K_{n} = k \giv K_{n+1} = k ) =  {2 (n - k + 1) \over  2 n - k + 1 }
\en
\eq
\lb{pred4444}
\PR ( K_{n} = k -1 \giv K_{n+1} = k ) =  {k - 1 \over  2 n - k + 1 } .
\en
As a check, the fact that $(K_n(\la), n \ge 1)$ has the same co-transition 
probabilities can be read from \re{pred33}, \re{pred44} and \re{fxy}.
It can be shown that the Markov chains $(K_n(\la), n \ge 1)$
for $\la \in [0,\infty]$, with definition by weak continuity for 
$\la = 0$ or $\infty$, are the extreme points of the convex set of all 
laws of Markov chains with these co-transition probabilities. A
generalization of this fact, to $\alpha \in (0,1)$ instead of 
$\alpha = \hf$, and similar considerations for $\alpha = 0$,
yield the second sentence of Theorem \ref{charc}.

To illustrate the formulas above, according to \re{cons3} and \re{momhf1},
or \re{eppf} for $n = 2$, given $L_1 = \ll$,
two independent uniform $[0,1]$ variables fall in the same excursion interval of the Brownian motion with probability
\eq
\lb{muone}
p_{\scrhf} ( 2 \bgiv \ll ) =  \mu(1 \bgiv \ll) = h_{-2} (\ll) = 1 - \ll h_{-1} (\ll)
\en
and in different excursion intervals with probability $\ll h_{-1} (\ll)$. 
According to \re{eppf} for $n=3$, given $L_1 = \ll$,
three independent uniform random points $U_1,U_2,U_3$ with uniform
distribution on $[0,1]$ fall 
in the same excursion interval of a Brownian motion or Brownian bridge
with probability
\eq
\lb{lxy}
\PR (K_3 (\ll) = 1 ) =
p_{\scrhf} (3 \bgiv \ll ) = 3 h_{-4} (\ll ) 
= 1 + \hf \ll ^2 - ( \thf \ll + \hf \ll ^3 ) h_{-1} (\ll )
\en
while $U_1$ and $U_2$ fall in one excursion interval and $U_3 $ in another
with probability
\eq
\lb{mxy}
\th \, \PR (K_3 (\ll) = 2 )  =
p_{\scrhf}(2,1 \bgiv \ll) = \ll h_{-3} (\ll) = - \hf \ll ^2 + ( \hf \ll + \hf \ll ^3 ) h_{-1}  (\ll )
\en
and the three points fall in three different excursion intervals with probability
\eq
\lb{nxy}
\PR (K_3 (\ll) = 3 ) =
p_{\scrhf} (1,1,1 \bgiv \ll ) =  \ll^2 h_{-2} (\ll) = \ll ^2 - \ll ^3  h_{-1}  (\ll ).
\en
As a check, the sum of expressions for $\PR (K_3 (\ll) = k )$ over
$k = 1,2,3$ reduces to 1.
Since
\eq
\lb{sumform}
\PR (K_n (\ll) = k  ) 
=
\sum_{n_1 \ge \cdots \ge n_k} 
\#(n_1, \cdots, n_k) p_{\scrhf}  (n_1, \cdots, n_k \bgiv \ll)  
\en
where the sum is over all decreasing sequences of positive
integers $(n_1, \cdots, n_k)$ with sum $n$, and
$\#(n_1, \cdots, n_k)$ is the number of distinct partitions of $\IN_n$ into
$k$ subsets of sizes $(n_1, \cdots, n_k)$,
formula \re{fxy} amounts to
\eq
\lb{hfiden}
\sum_{n_1 \ge \cdots \ge n_k} 
\#(n_1, \cdots, n_k) 
\prod_{i= 1}^k [\hf ]_{n_i - 1} = 
{ 2n - k - 1 \choose n - 1 } \, { \Gamma (n ) \over \Gamma (k ) } \, 2 ^{ 2 k - 2n } 
\en
which can be checked as follows.
According to \re{class.esf} and \re{pkhf1}, the unconditional EPPF of the Brownian excursion partition $\Pi$ is
\eq
\lb{brownepf}
p_{\scrhf,0}(n_1, \cdots, n_k) = { \Gamma ( k ) \over 2^{k-1} \Gamma (n ) }
\prod_{i= 1}^k [\hf ]_{n_i - 1}
\en
so \re{hfiden} can be deduced from \re{brownepf}, \re{brownkn}, and the unconditioned form of \re{sumform}.

\subsection{Some identitities}
As a consequence of \re{pdint} and \re{momhf1}, for all $q > - \hf$ and $\theta > - \hf$ there is the identity
\eq
\lb{strq}
{2 \over \ER ( |B_1|^{2 \theta } ) } \int_0^\infty \ll^{2 \theta} \mu( q \bgiv \ll ) \phi(\ll) d \ll
= { \Gamma (\theta + 1) \Gamma ( q + \hf ) \over \Gamma(\hf) \Gamma( q + \theta + 1)} 
\en
where the right side is the $q$th moment of the beta$(\hf, \hf + \theta)$ structural distribution of
$\PD(\hf,\theta)$, and on the left side this moment is computed by conditioning on $L_1$.
As in \re{betatf}, for each fixed $q$ this formula provides a Mellin transform which
uniquely determines $\mu( q \bgiv \ll )$ as a function of $\ll$.
In view of \re{strq} and \re{bmoms},
the formula \re{momhf1} for $\mu( q \bgiv \ll )$ in terms of 
the Hermite function amounts to the identity
\eq
\lb{strq3}
2 \int_0^\infty \ll^{2 \theta} h_{- 2 q } (\ll) \phi(\ll) d \ll
=  2 ^{ -\theta - q } { \Gamma( 2 \theta + 1)  \over \Gamma (q + \theta + 1) } .
\en
As checks, since $h_0(x) = 1$ and $h_{-1}(x) = \Phi(x)/\varphi(x)$,
the case $q = 0$ is obvious,
and the case $q = \hf$ is easily verified since then the left side of 
\re{strq} equals $(2 \theta + 1)^{-1} \ER ( |B_1|^{2 \theta + 1} )$ by integration by parts. Formula \re{strq3} can then be verified for $q = m/2$ for all 
$m = 0,1,2, \ldots$, using the recursion \re{recurh}.
Formula \re{strq3} was just derived for $q > -\hf$, but both sides are entire
functions of $q$, so the identity holds for all $q \in \complex$.
Using the series formula \re{hermq} and integrating term by term, the substitution $r = \theta + \hf$ allows the
identity \re{strq3} to be rewritten in the symmetric form
\eq
\lb{idsymm}
\sum_{j = 0}^\infty  \Gamma \left( q + {j \over 2} \right) \Gamma \left(r + {j \over 2} \right) { (-2)^j \over j! } = 
{4 \sqrt{ \pi} \Gamma( 2 q) \Gamma( 2 r ) \over \Gamma ( q + r + 1/2) } 
\en
where the series is absolutely convergent for real $q$ and $r$ with $q + r + \hf < -1$, and can otherwise
be summed by Abel's method provided neither $2q$ nor $2r$ is a non-positive integer. This version of the identity
is easily verified using standard identities involving Gauss's hypergeometric function and the gamma function.
For $-2q = n$ a positive integer, when $h_n$ is the $n$th Hermite polynomial
$$
h_n(x) = \sum_{k = 0}^{ \lfloor n/2 \rfloor } h_{n,k} \, x^{n - 2k } \mbox{ with } h_{n,k} = (-1)^k { n  \choose 2 k } { ( 2 k )! \over 2^{k} k! } 
$$
the identity \re{strq3} reduces easily to the following pair of identities
of polynomials in $\theta$, 
which relate the rising and falling factorials 
$[x]_n:= x(x+1) \cdots (x + n-1)$ and $(x)_n:= x(x-1) \cdots (x - n+1)$,
and which are easily verified directly: for $m = 0,1, 2 \ldots$
$$
\sum_{k = 0}^m h_{2m,k} 2^{-k} [ \theta + \hf]_{m-k} = (\theta)_m
$$
and 
$$
\sum_{k = 0}^m h_{2m+1,k} 2^{-k} [ \theta + 1]_{m-k} = (\theta - \hf)_m .
$$
Thus the coefficients of the Hermite polynomials are related to some instances of generalized Stirling numbers \cite{hsu98,csp}.

\vspace{.25 in}
\noindent
{\Large \bf Acknowledgment}

\noindent
Thanks to Gr\'egory Miermont for careful reading of a draft of this paper.

\def\cprime{$'$} \def\polhk#1{\setbox0=\hbox{#1}{\ooalign{\hidewidth
  \lower1.5ex\hbox{`}\hidewidth\crcr\unhbox0}}} \def\cprime{$'$}
  \def\cprime{$'$} \def\cprime{$'$}
  \def\polhk#1{\setbox0=\hbox{#1}{\ooalign{\hidewidth
  \lower1.5ex\hbox{`}\hidewidth\crcr\unhbox0}}} \def\cprime{$'$}
  \def\cprime{$'$} \def\polhk#1{\setbox0=\hbox{#1}{\ooalign{\hidewidth
  \lower1.5ex\hbox{`}\hidewidth\crcr\unhbox0}}} \def\cprime{$'$}
  \def\cprime{$'$} \def\cydot{\leavevmode\raise.4ex\hbox{.}} \def\cprime{$'$}
  \def\cprime{$'$} \def\cprime{$'$} \def\cprime{$'$}

%
%

\affiliation{Jim Pitman, Department of Statistics, University of California, Berkeley \\ {\tt pitman@stat.berkeley.edu}}



\begin{thebibliography}{10}

\bibitem{as65}
M.~Abramowitz and I.~A. Stegun.
\newblock {\em Handbook of Mathematical Functions}.
\newblock Dover, New York, 1965.

\bibitem{ap92}
D.~Aldous and J.~Pitman.
\newblock Brownian bridge asymptotics for random mappings.
\newblock {\em Random Structures and Algorithms}, 5:487--512, 1994.

\bibitem{jpda97sac}
D.J. Aldous and J.~Pitman.
\newblock The standard additive coalescent.
\newblock {\em Ann. Probab.}, 26:1703--1726, 1998.

\bibitem{abt01}
R.~A. Arratia, A.~D. Barbour, and S.~Tavar{\'e}.
\newblock {Logarithmic combinatorial structures: a probabilistic approach}.
\newblock Book in preparation. Available via {\tt
  www-hto.usc.edu/books/tavare/ABT/index.html}, 2001.

\bibitem{bertoin00f}
J.~Bertoin.
\newblock A fragmentation process connected to {B}rownian motion.
\newblock {\em Probab. Theory Related Fields}, 117(2):289--301, 2000.

\bibitem{bert02ss}
J.~Bertoin.
\newblock Self-similar fragmentations.
\newblock {\em Ann. Inst. H. Poincar\'e Probab. Statist.}, 38(3):319--340,
  2002.

\bibitem{bl73}
D.~Blackwell and J.B. MacQueen.
\newblock {Ferguson distributions via P\'{o}lya urn schemes}.
\newblock {\em Ann. Statist.}, 1:353--355, 1973.

\bibitem{brix99}
A.~Brix.
\newblock Generalized gamma measures and shot-noise {C}ox processes.
\newblock {\em Adv. in Appl. Probab.}, 31(4):929--953, 1999.

\bibitem{chassjan01}
P.~Chassaing and S.~Janson.
\newblock A {V}ervaat-like path transformation for the reflected {B}rownian
  bridge conditioned on its local time at 0.
\newblock {\em Ann. Probab.}, 29(4):1755--1779, 2001.

\bibitem{chasslou99}
P.~Chassaing and G.~Louchard.
\newblock Phase transition for parking blocks, {B}rownian excursion and
  coalescence.
\newblock {\em Random Structures Algorithms}, 21:76--119, 2002.

\bibitem{derrida81}
B.~Derrida.
\newblock Random-energy model: an exactly solvable model of disordered systems.
\newblock {\em Phys. Rev. B (3)}, 24(5):2613--2626, 1981.

\bibitem{derrida94}
B.~Derrida.
\newblock Non-self-averaging effects in sums of random variables, spin glasses,
  random maps and random walks.
\newblock In M.~Fannes, C.~Maes, and A.~Verbeure, editors, {\em On Three
  Levels: Micro-, Meso-, and Macro-Approaches in Physics}, NATO ASI Series,
  pages 125--137. Plenum Press, New York and London, 1994.

\bibitem{derrida97}
B.~Derrida.
\newblock From random walks to spin glasses.
\newblock {\em Phys. D}, 107(2-4):186--198, 1997.
\newblock Landscape paradigms in physics and biology (Los Alamos, NM, 1996).

\bibitem{doetsch37}
G.~Doetsch.
\newblock {\em {Theorie und Anwendung der Laplace-Transformation}}.
\newblock Berlin, 1937.

\bibitem{en78}
S.~Engen.
\newblock {\em Stochastic Abundance Models with Emphasis on Biological
  Communities and Species Diversity}.
\newblock Chapman and Hall Ltd., 1978.

\bibitem{erdelyi53II}
A.~Erd{\'e}lyi et~al.
\newblock {\em {Higher Transcendental Functions}}, volume~II of {\em Bateman
  Manuscript Project}.
\newblock McGraw-Hill, New York, 1953.

\bibitem{et95}
W.~J. Ewens and S.~Tavar{\'e}.
\newblock {The Ewens sampling formula}.
\newblock In N.~S. Johnson, S.~Kotz, and N.~Balakrishnan, editors, {\em
  Multivariate Discrete Distributions}. Wiley, New York, 1995.

\bibitem{ew72}
W.J. Ewens.
\newblock The sampling theory of selectively neutral alleles.
\newblock {\em Theor. Popul. Biol.}, 3:87 -- 112, 1972.

\bibitem{grotespeed02}
M.~Grote and T.~P. Speed.
\newblock {Approximate Ewens formulae for symmetric overdominance selection}.
\newblock {\em Annals of Applied Probability}, 12: 637 -- 663, 2002.

\bibitem{hansenp98}
B.~Hansen and J.~Pitman.
\newblock {Prediction rules and exchangeable sequences related to species
  sampling}.
\newblock {\em Stat. and Prob. Letters}, 46(520):251--256, 2000.

\bibitem{ho87}
F.~M. Hoppe.
\newblock The sampling theory of neutral alleles and an urn model in population
  genetics.
\newblock {\em Journal of Mathematical Biology}, 25:123 -- 159, 1987.

\bibitem{hsu98}
L.~C. Hsu and P.~J.-S. Shiue.
\newblock A unified approach to generalized {S}tirling numbers.
\newblock {\em Adv. in Appl. Math.}, 20(3):366--384, 1998.

\bibitem{james02pp}
L.~F. James.
\newblock {Poisson process partition calculus with applications to exchangeable
  models and Bayesian nonparametrics}.
\newblock arXiv:math.PR/0205093, 2002.

\bibitem{james02}
L.F. James.
\newblock Bayesian calculus for gamma processes with applications to
  semiparametric intensity models.
\newblock {\em Sankhy{\=a}}, 2002.
\newblock to appear.

\bibitem{janson01h}
S.~Janson.
\newblock Asymptotic distribution for the cost of linear probing hashing.
\newblock {\em Random Structures Algorithms}, 19(3-4):438--471, 2001.
\newblock Analysis of algorithms (Krynica Morska, 2000).

\bibitem{JK70}
N.~L. Johnson and S.~Kotz.
\newblock {\em Continuous Univariate Distributions, volume 2}.
\newblock Wiley, 1970.

\bibitem{kar67urn}
S.~Karlin.
\newblock Central limit theorems for certain infinite urn schemes.
\newblock {\em J. Math. Mech.}, 17:373--401, 1967.

\bibitem{kerov95}
S.~Kerov.
\newblock {Coherent random allocations and the Ewens-Pitman formula}.
\newblock PDMI Preprint, Steklov Math. Institute, St. Petersburg, 1995.

\bibitem{ki75}
J.~F.~C. Kingman.
\newblock Random discrete distributions.
\newblock {\em J. Roy. Statist. Soc. B}, 37:1--22, 1975.

\bibitem{ki78b}
J.~F.~C. Kingman.
\newblock The representation of partition structures.
\newblock {\em J. London Math. Soc.}, 18:374--380, 1978.

\bibitem{ki82co}
J.~F.~C. Kingman.
\newblock The coalescent.
\newblock {\em Stochastic Processes and their Applications}, 13:235--248, 1982.

\bibitem{king82}
J.~F.~C. Kingman.
\newblock Exchangeability and the evolution of large populations.
\newblock In G.~Koch and F.~Spizzichino, editors, {\em Exchangeability in
  probability and statistics (Rome, 1981)}, pages 97--112. North-Holland,
  Amsterdam, 1982.

\bibitem{ki93}
J.~F.~C. Kingman.
\newblock {\em Poisson Processes}.
\newblock Clarendon Press, Oxford, 1993.

\bibitem{ksor97}
U.~K{\"u}chler and M.~S{\o}rensen.
\newblock {\em Exponential families of stochastic processes}.
\newblock Springer-Verlag, New York, 1997.

\bibitem{lebedev65}
N.~N. Lebedev.
\newblock {\em Special Functions and their Applications}.
\newblock Prentice-Hall, Englewood Cliffs, N.J., 1965.

\bibitem{lev39}
P.~L{\'e}vy.
\newblock Sur certains processus stochastiques homog{\`e}nes.
\newblock {\em Compositio Math.}, 7:283--339, 1939.

\bibitem{mc65}
J.~W. McCloskey.
\newblock A model for the distribution of individuals by species in an
  environment.
\newblock Ph. D. thesis, Michigan State University, 1965.

\bibitem{miersch01}
G.~Miermont and J.~Schweinsberg.
\newblock Self-similar fragmentations and stable subordinators.
\newblock Technical Report 726, 2001.
\newblock Pr\'epublication du Laboratoire de Probabilit\'es et mod\`eles
  al\'eatoires, Universite Paris VI. Available via {\tt
  http://www.proba.jussieu.fr}.

\bibitem{miller55}
J.~C.~P. Miller.
\newblock {\em {Tables of Weber parabolic cylinder functions}}.
\newblock Her Majesty's Stationery Office, London, 1955.

\bibitem{per91}
M.~Perman.
\newblock Order statistics for jumps of normalized subordinators.
\newblock {\em Stoch. Proc. Appl.}, 46:{267--281}, 1993.

\bibitem{ppy92}
M.~Perman, J.~Pitman, and M.~Yor.
\newblock Size-biased sampling of {P}oisson point processes and excursions.
\newblock {\em Probab. Th. Rel. Fields}, 92:21--39, 1992.

\bibitem{pitman95pk}
J.~Pitman.
\newblock {Poisson-Kingman partitions}.
\newblock Preprint available via {\tt www.stat.berkeley.edu/users/pitman},
  1995.

\bibitem{jp.epe}
J.~Pitman.
\newblock Exchangeable and partially exchangeable random partitions.
\newblock {\em Probab. Th. Rel. Fields}, 102:145--158, 1995.

\bibitem{jp.isbp}
J.~Pitman.
\newblock Random discrete distributions invariant under size-biased
  permutation.
\newblock {\em Adv. Appl. Prob.}, 28:525--539, 1996.

\bibitem{jp96bl}
J.~Pitman.
\newblock {Some developments of the Blackwell-MacQueen urn scheme}.
\newblock In T.S.~Ferguson et~al., editor, {\em {Statistics, Probability and
  Game Theory; Papers in honor of David Blackwell}}, volume~30 of {\em Lecture
  Notes-Monograph Series}, pages 245--267. Institute of Mathematical
  Statistics, Hayward, California, 1996.

\bibitem{jp.bmpart}
J.~Pitman.
\newblock {Partition structures derived from Brownian motion and stable
  subordinators}.
\newblock {\em Bernoulli}, 3:79--96, 1997.

\bibitem{jp99sam}
J.~Pitman.
\newblock {Brownian motion, bridge, excursion and meander characterized by
  sampling at independent uniform times}.
\newblock {\em Electron. J. Probab.}, 4:Paper 11, 1--33, 1999.

\bibitem{csp}
J.~Pitman.
\newblock {Combinatorial Stochastic Processes}.
\newblock Technical Report 621, Dept. Statistics, U.C. Berkeley, 2002.
\newblock Lecture notes for St. Flour course, July 2002. Available via
  {\tt www.stat.berkeley.edu}.

\bibitem{py92}
J.~Pitman and M.~Yor.
\newblock Arcsine laws and interval partitions derived from a stable
  subordinator.
\newblock {\em Proc. London Math. Soc. (3)}, 65:326--356, 1992.

\bibitem{py95pd2}
J.~Pitman and M.~Yor.
\newblock {The two-parameter Poisson-Dirichlet distribution derived from a
  stable subordinator}.
\newblock {\em Ann. Probab.}, 25:855--900, 1997.

\bibitem{py02h}
J.~Pitman and M.~Yor.
\newblock {Brownian interpretations of Mill's ratio, Hermite and parabolic
  cylinder functions}.
\newblock In preparation, 2002.

\bibitem{po46}
H.~Pollard.
\newblock {The representation of $e^{-x^\la}$ as a Laplace integral}.
\newblock {\em Bull. Amer. Math. Soc.}, 52:908--910, 1946.

\bibitem{ruelle87}
D.~Ruelle.
\newblock {A mathematical reformulation of Derrida's REM and GREM}.
\newblock {\em Comm. Math. Phys.}, 108:225--239, 1987.

\bibitem{seshadri93}
V.~Seshadri.
\newblock {\em The inverse {G}aussian distribution}.
\newblock The Clarendon Press, Oxford University Press, New York, 1993.

\bibitem{tala01}
M.~Talagrand.
\newblock {\em {Spin glasses: a challenge to mathematicians}}.
\newblock Springer, 2001.
\newblock Book in preparation.

\bibitem{toscano71p}
L.~Toscano.
\newblock Sulle funzioni del cilindro parabolico.
\newblock {\em Matematiche (Catania)}, 26:104--126 (1972), 1971.

\bibitem{za94}
S.L. Zabell.
\newblock The continuum of inductive methods revisited.
\newblock In J.~Earman and J.~D. Norton, editors, {\em The Cosmos of Science},
  Pittsburgh-Konstanz Series in the Philosophy and History of Science, pages
  351--385. University of Pittsburgh Press/Universit{\"a}tsverlag Konstanz,
  1997.

\end{thebibliography}
\end{document}